\documentclass[11pt,a4paper]{article}
\usepackage{color,pdfsync}
\usepackage{times,amsmath,amssymb,amsfonts,lscape,multirow}
\usepackage{url}
\usepackage{a4wide}
\usepackage{epsfig,epic,eepic,graphicx,latexsym,pdfsync}
\newtheorem{theorem}{Theorem}
\newtheorem{lemma}[theorem]{Lemma}
\newtheorem{proposition}[theorem]{Proposition}
\newtheorem{corollary}[theorem]{Corollary}

\newcommand{\alert}{}

\newtheorem{definition}{Definition}

\makeatletter

\@addtoreset{equation}{section}
\makeatother

\newcommand{\N}{{\mathbb N}}

\newcommand{\R}{{\mathbb R}}

\newcommand{\pa}{{\partial}}
\newcommand{\na}{{\nabla}}

\newcommand{\eps}{{\varepsilon}}

\def\div{\hbox{div  }}

\title{Existence of weak solutions up to collision \\
 for viscous fluid-solid systems with slip}
\author{David G\'erard-Varet\footnote{Institut de Math\'ematiques de Jussieu, 175 rue du Chevaleret, 75013 Paris. Mail: \texttt{gerard-varet@math.jussieu.fr}} \hspace{0.5cm} Matthieu Hillairet\footnote{Ceremade, Place du Mar\'echal de Lattre de Tassigny, 75775 Paris Cedex 16. Mail: \texttt{hillairet@ceremade.dauphine.fr}}}

\begin{document}
\maketitle
\begin{abstract}
We study in this paper the movement of a rigid solid inside an incompressible Navier-Stokes flow, within a bounded domain. We consider the case where slip is allowed at the fluid/solid interface, through a Navier condition. Taking into account slip at the interface is very natural within this model, as  classical no-slip conditions lead to  unrealistic collisional behavior  between the solid and the domain boundary. We prove for this model existence of weak solutions of Leray type,  up to collision, in three dimensions. The key point is that, due to the slip condition,  the velocity field is discontinuous across the fluid/solid interface. This prevents from obtaining global $H^1$ bounds on the velocity, which makes many aspects of the theory of weak solutions for Dirichlet conditions inappropriate.  
\end{abstract}

\section{Introduction}
The general concern of this paper is the dynamics of solid bodies in a fluid flow. This dynamics is relevant to many  natural and industrial processes, like  blood flows, sprays, or  design of micro swimmers. 

\medskip
A main problem to understand this dynamics is to compute the drag exerted by the flow on the bodies. From the mathematical point of view, a natural approach to this problem is  to use the Euler or Navier-Stokes equations to model the flow. However, this generates serious difficulties. A famous one is {\em D'Alembert's paradox}, related to the Euler equation: in an incompressible and inviscid potential flow, a solid body undergoes no drag \cite{MaPu}.

\medskip
But  the Navier-Stokes equations also raise modeling issues. Let us consider for instance a single solid moving in a viscous fluid. We denote by 
 $S(t) \subset \R^3$, $\: F(t) \subset \R^3$ the solid and fluid domains at time $t$, and $\displaystyle \Omega \: := \: \overline{S(t)} \, \cup \, F(t)$
   the total domain. We assume that the fluid is governed by the Navier-Stokes equations.  We denote  $u_F$ and  $p_F$ its velocity and internal pressure, $\rho_F$ its density, $\mu_F$ its  viscosity. Thus:
\begin{equation} \label{NS} 
\left\{
\begin{aligned}
 \rho_F \left( \pa_t u_F + u_F \cdot \na u_F \right) - \mu_F \Delta u_F
 \: & =
\:  - \na p_F - \rho_F g, \quad t > 0, \: x \in F(t), \\ 
 \div u_F \: &  = \:  0, \hspace{1cm} t > 0, \: x \in F(t), 
\end{aligned}
\right.
\end{equation}  
 with  $- \rho_F g$ the gravitational force. In parallel to the fluid modeling, we  write the conservation of linear and angular  momentum for the body. Denoting $x_S(t)\in \R^3$ the position of the center of mass, $\: U_S(t) \in \R^3$ its velocity,  and 
  $\: \omega_S(t) \in \R^3$ the angular velocity at time  $t$, these conservation laws read
 \begin{equation} \label{S3D} 
\left\{
\begin{aligned} m_S \, \frac{d}{dt}U_S(t)  & \: = \: - \int_{\pa S(t)} \Sigma_F \nu \, d\sigma
 \: - \: m_S \, g, \\ 
\frac{d}{dt} \left( J_S \, \omega_S(t) \right) & \: = \: - \:    \int_{\pa S(t)} \left( x - x_S(t)
 \right) \times (\Sigma_F \, \nu) \, d\sigma  \:  + \: \rho_S  \int_{S(t)}   \left( x - x_S(t) \right)  \times  (-g \, ).   
\end{aligned}
\right.
\end{equation}
Following  standard notations,   $\rho_S$ and $m_S \, :=  \, \rho_S  |S(0)|$ are the density and mass of the solid (independent of $t$ and $x$),    
$\Sigma_F(t,x) \in M_3(\R)$ is the newtonian tensor of the fluid:
$$  \Sigma_F    \:  = \:  \left( 2 \mu_F D(u_F)  - p_F \, I_d \, \right), $$
and $J_S(t)  \in M_3(\R) $ is the inertia matrix of the solid: 
$$J_S(t) \: := \: \rho_S \, \int_{S(t)}\left(  |x - x_S(t)|^2 \, I_d  \, - \,   (x - x_S(t) ) \otimes (x - x_S(t) ) \right)  dx.$$
 The vector  $\nu=\nu(t, \cdot)$ is the unit normal vector pointing inside the solid $S(t)$. Note that the 
 velocity  $u_S(t,x)$ at each point $x$ of the solid reads 
$$ u_S(t,x) \:  := \:  U_S(t) \: + \: \omega_S(t) \times (x - x_S(t)). $$

\medskip
To close the system,  one usually imposes  no-slip conditions, both  at the fluid-solid interface and the cavity boundary: 
\begin{equation} \label{noslip} 
\left\{
\begin{aligned}
& u_F\vert_{\pa S(t)} \:  = \: u_S\vert_{\pa S(t)} \\
& u_F\vert_{\pa \Omega} =  0, 
\end{aligned}
\right.
 \end{equation}
and one specifies  the initial data:  the initial position of the solid $S_0$, 
$$ \: u_{F,0} \: :=  \: u_F\vert_{t=0} \:  \mbox{ and } \: 
u_{S,0} \: := \: U_{S,0} +  \omega_S(0) \times (x- x_{S_0}).$$ 

\medskip
One could believe that system \eqref{NS}-\eqref{S3D}-\eqref{noslip} is a good model for the interaction between a solid and a viscous fluid. Far from it:  in the case of a sphere falling over a flat wall
$$ S(0) \: :=  \: e_3 + B(0,1/2), \quad \Omega \: := \: \{ x_3 > 0 \}, $$
it predicts that no collision is possible between the solid and the wall !  This no-collision paradox has been known from specialists since the 1960's, after articles by Cox and Brenner \cite{Brenner&Cox63} and Cooley and O'Neill \cite{Cooley&Oneill69} in the context of Stokes equations.  Since then, the no-collision paradox has been  confirmed at the level of the Navier-Stokes equations  (see  \cite{Hillairet07,HiTa}, and the preliminary result in \cite{SanMartin&al02}).  

\medskip
Of course, such a  result is unrealistic, as it goes against Archimedes' principle. It suggests that the Navier-Stokes equations are not relevant to collisional and post-collisional descriptions. Hence, many physicists have tried to find an explanation for the paradox.  {\em We shall focus here on one possible explanation, namely the no-slip condition}. The idea is that, when the distance between the solids gets very small (below the micrometer), the no-slip condition is no longer accurate, and must be replaced by a Navier condition:  
 \begin{equation} \label{slip}
\left\{
\begin{aligned}
 \bigl(u_F -  \, u_S \bigr) \cdot \nu\vert_{\pa S(t)}  & = \:  0,  &  \quad   \bigl(u_F -  \, u_S\bigr)  \times \nu\vert_{\pa S(t)} & = \:   \,    - 2 \beta_S (\Sigma_F \nu)    \times \nu\vert_{\pa S(t)}, \\
u_F  \cdot \nu\vert_{\pa \Omega}  & = \:  0,  & \quad u_F  \times \nu\vert_{\pa \Omega}  & = \:  -2 \beta_\Omega \,   (\Sigma_F \nu)    \times \nu\vert_{\pa \Omega}.
\end{aligned} 
\right.
\end{equation}
In other words, only the normal component of the relative velocity of the fluid is zero, to ensure impermeability. The tangential ones are non-zero, and proportional to the stress constraint, with constant slip lengths $\beta_S, \beta_\Omega > 0$. For a recent discussion of the Navier condition, notably in the context of microfluidics, we refer to \cite{Lauga:2007}. See also the seminal paper \cite{Hocking}. Let us point out that the Navier-condition is sometimes used as  a wall law, to describe the averaged effect of rough hydrophobic surfaces \cite{Bocquet:2007}. 

\medskip
The effect of slip conditions \eqref{slip} on collision was investigated recently  by the authors in article \cite{GeHi2}. 
More precisely, this article is devoted to a simplified model, in which 
\begin{itemize}
\item The Navier-Stokes equations are replaced by the steady Stokes ones (quasi-static regime). 
\item The domain $\Omega$ is a half-space, the solid $S$ is a sphere. 
\end{itemize}
In this context, denoting $h(t)$ the distance between $S(t)$ and the plane $\pa \Omega$, it is shown that the dynamics obeys the reduced ODE 
$$  \ddot{h} \: = \: \dot{h} \, {\cal D}(h) \: + \:  \frac{(\rho_F - \rho_S)}{\rho_S} \,  g$$
where the drag term ${\cal D}(h)$ satisfies ${\cal D}(h) = O(|\ln h|)$ as $ h \rightarrow 0$. This is in sharp contrast with the case of no-slip conditions \eqref{noslip}, for which $ {\cal D}(h) \sim \frac{C}{h}$. In particular, it allows for collisions in finite time. We refer to \cite{GeHi2} for all details and other results in the context of rough boundaries.  

\medskip
Hence, paper \cite{GeHi2} provides a resolution of the paradox: one can {\it a priori} keep the Navier-Stokes equations, up to considering the Navier boundary conditions \eqref{slip}.   Nevertheless, the analysis in \cite{GeHi2} is limited to simple configurations and to  Stokes flows. In the context of the full 3D Navier-Stokes system \eqref{NS}, more complicated behaviors may occur. For instance, smooth solutions may exhibit singularities prior to  any collision. To describe the qualitative features of the collision,  {\em one needs to consider weak solutions}. The theory of weak solutions is well understood in the case of no-slip conditions and many references will be given in the next section. {\em However, the existence of weak solutions with Navier conditions  has been so far an open question, due to serious additional mathematical difficulties. To address this question is the purpose of the present paper}. Broadly, we shall build weak solutions for system \eqref{NS}-\eqref{S3D}-\eqref{slip}, up to collision between the solid and the cavity $\Omega$. 

\medskip
The paper is organized as follows. Section \ref{results} contains the statement of our main result: we give a definition of weak solutions, and state the existence of such solutions as long as no contact occurs. We explain the main difficulties in proving their existence, in comparison to the results available for no-slip conditions. We conclude Section \ref{results} by an outline of our proof, to be  carried out in sections \ref{transport} to \ref{convergence}. More precisely: 
\begin{itemize}
\item Section \ref{transport} is devoted to an auxiliary nonlinear transport equation, which is crucial to our approximation procedure. 
\item Section \ref{approx} is dedicated to the construction of solutions for well-chosen approximations of the Navier-Stokes / solid dynamics. 
\item Section \ref{convergence} describes the limit procedure, from the approximate to the exact system.
\end{itemize}

\section{Main result and ideas} \label{results}

\subsection{Weak solutions of Navier-Stokes with slip conditions}
The aim of this paragraph is to define a weak formulation and weak solutions for system  \eqref{NS}-\eqref{S3D}-\eqref{slip}, that is in the case of slip conditions of Navier type. 
We remind that in the case of no-slip conditions, the theory of weak solutions has been successfully achieved over the last ten years, first up to collision (see \cite{Desjardins&Esteban99})  and then globally in time (see \cite{SanMartin&al02} in the 2D case,  \cite{Feireisl03} in the 3D case). Let us also mention the alternative approach in  \cite{BostCottetMaitre10}, and the recent result \cite{GlassSueurpp} on the uniqueness of 2D weak solutions up to collision. 

\medskip
 As usual,  in order to derive a weak formulation, the starting point  is formal multiplication by appropriate test functions. These test functions must of course look like the solution itself. Notably, they must be rigid vector fields in the solid domain $S$. This forces the space of test functions to depend on the solution itself: it is a classical difficulty, already recognized in the no-slip case.  {\em A key feature of the slip conditions is that these test functions, and also the solution, will be moreover discontinuous across the fluid/solid interface. } Indeed, the first line of \eqref{slip} ensures the continuity of the normal component, but the tangential ones may have a jump. It is a strong difference with regards to boundary conditions  \eqref{noslip}, and it  will generate many difficulties throughout the paper. 
 
 \medskip
 We first introduce some notation for the  classical spaces of solenoidal vector fields.  Let ${\cal O}$ be a  Lipschitz domain. We set 
 \begin{equation*}
\begin{aligned}
& {\cal D}_\sigma({\cal O}) \: := \: \left\{ \varphi \in  {\cal D}({\cal O}), \quad \div \varphi = 0  \right\}, \quad 
    {\cal D}_\sigma(\overline{{\cal O}}) \: := \: \left\{ \varphi\vert_{{\cal O}}, \: \varphi  \in  {\cal D}_\sigma(\R^3)  \right\},
 \\
& L^2_\sigma({\cal O}) \: := \: \:  \mbox{  the closure of    ${\cal D}_\sigma({\cal O})$ in $L^2({\cal O})$}, \quad  H^1_\sigma({\cal O}) \: := \: H^1({\cal O}) \cap L^2_\sigma({\cal O}), \\
& H^1_\sigma(\overline{{\cal O}})  \: := \:  \:  \mbox{  the closure of    ${\cal D}_\sigma(\overline{{\cal O}})$ in $H^1({\cal O})$}
 \end{aligned}
 \end{equation*} 
We remind that elements $u$ of $L^2_\sigma({\cal O})$ satisfy $u \cdot \nu = 0$ in $H^{-1/2}({\pa \cal O})$. 

\medskip
 We also introduce the finite dimensional space of rigid vector fields in $\R^3$: 
 $${\cal R} \: := \: \{\varphi_s, \quad  \varphi_s(x) = V + \omega \times x, \quad \mbox{ for some } \: V \in \R^3, \:  \omega \in \R^3 \}. $$
Finally, we  define for  any $T > 0$  the space ${\cal T}_T$ of test functions over $[0,T)$: 
\begin{multline*}
{\cal T}_T \: := \: \Bigl\{ \varphi \in C([0,T]; L^2_\sigma(\Omega)), \: \mbox{there exists }  \: \varphi_F \in {\cal D}([0,T);  {\cal D}_\sigma(\overline{\Omega})), \:  \varphi_S \in {\cal D}([0,T);  {\cal R})  \\
 \mbox{ such that } \:  \varphi(t,\cdot) = \varphi_F(t,\cdot) \:  \mbox{ on } F(t), \:   \varphi(t,\cdot) = \varphi_S(t,\cdot) \:  \mbox{ on } S(t), \: \mbox{ for all } \:  t \in [0,T] \Bigr\}.
\end{multline*}
Let us point out once again that this space of test functions depends on the solution itself through the domains $S(t)$ and $F(t)$. Let us also notice that the constraint $\varphi(t, \cdot) \in L^2_\sigma(\Omega)$ encodes in a weak form the continuity of the normal component at $\pa S(t)$: 
$$ \varphi_F(t,\cdot) \cdot \nu = \varphi_S(t,\cdot) \cdot \nu \: \mbox{ at }\: \pa S(t), \quad \forall t \in [0,T). $$

Multiplying  \eqref{NS} by $\varphi \in {\cal T}_T$, integrating over $F(t)$, and integrating by parts,  we obtain (formally)
\begin{multline*}
\frac{d}{dt}\int_{F(t)}\rho_F  \, u_F \cdot \varphi_F -   \int_{F(t)} \rho_F \, u_F \cdot \pa_t \varphi_F - \int_{F(t)} \rho_F \, u_F \otimes u_F : \na \varphi_F \: + \: \int_{F(t)} \hspace{-0.2cm} 2 \mu_F D(u_F) : D(\varphi_F)  \\
\: = \:   \int_{\pa \Omega} \left(\Sigma_F \nu  \right)\cdot \varphi_F \: + \: 
\int_{\pa S(t)}   \left( \Sigma_F \nu \right) \cdot \varphi_F \: + \:    \int_{F(t)} \rho_F(-g) \cdot \varphi_F,  
\end{multline*}
where the normal vectors $\nu$, in the integrals at the right-hand side,   point resp.  outside $\Omega$ and inside  $S(t)$. Using \eqref{slip}:
\begin{align*}
\int_{\pa \Omega} \left(\Sigma_F \nu \right)\cdot \varphi_F & = -\frac{1}{2\beta_\Omega} \int_{\pa \Omega}  (u_F \times \nu) \cdot (\varphi_F \times \nu), \\ 
\int_{\pa S(t)}   \left( \Sigma_F \nu\right) \cdot \varphi_F & = -\frac{1}{2\beta_S} \int_{\pa S(t)}  ((u_F - u_S) \times \nu) \cdot ((\varphi_F - \varphi_S) \times \nu) \: + \: \int_{\pa S(t)}   \left(\Sigma_F \nu \right) \cdot \varphi_S 
 \end{align*}
Eventually, one can use \eqref{S3D} to write differently the last integral: tedious but straightforward  manipulations yield 
$$ \int_{\pa S(t)}  \left(\Sigma_F \nu \right) \cdot \varphi_S  =   -\frac{d}{dt} \int_{S(t)}\rho_S  \, u_S \cdot \varphi_S  + \int_{S(t)} \rho_S \, u_S \cdot \pa_t \varphi_S + \int_{S(t)} \rho_S (-g) \cdot \varphi_S. $$
Combining the previous identities and integrating from $0$ to $T$ entails 
\begin{equation} \label{momentum}
\begin{aligned}
& -   \int_0^T \int_{F(t)} \rho_F \, u_F \cdot \pa_t \varphi_F - \int_0^T \int_{S(t)}  \rho_S \, u_S \cdot \pa_t \varphi_S +  \int_0^T \int_{F(t)} \rho_F \, u_F \otimes u_F : \na \varphi_F \\ 
&  + \: \int_0^T \int_{F(t)} \hspace{-0.2cm} 2 \mu_F D(u_F) : D(\varphi_F) 
+ \frac{1}{2\beta_\Omega}\int_0^T \int_{\pa \Omega}  (u_F \times \nu) \cdot (\varphi_F \times \nu)
\\
& + \frac{1}{2\beta_S} \int_0^T  \int_{\pa S(t)}  ((u_F - u_S) \times \nu) \cdot ((\varphi_F - \varphi_S) \times \nu) 
\\
& = \int_0^T  \int_{F(t)} \rho_F(-g) \cdot \varphi_F +  \int_0^T \int_{S(t)} \rho_S (-g) \cdot \varphi_S \\
& + \int_{F(0)} \rho_F u_{F,0} \cdot \varphi_F\vert_{t=0} +  \int_{S(0)} \rho_S u_{S,0} \cdot \varphi_S\vert_{t=0} 
  \end{aligned}
  \end{equation}
Equation \eqref{momentum} is a global weak formulation of the momentum equations \eqref{NS} and \eqref{S3D}, taking the slip conditions \eqref{slip} into account. 
Setting $\varphi = u$ in the above formal computations yields that, for all $t \in [0,T] :$
\begin{equation} \label{mainenergybound}
 \begin{aligned} & \int_{F(t)} \frac{1}{2}\rho_F |u_F(t, \cdot)|^2 + \int_{S(t)}\frac{1}{2} \rho_S |u_S(t, \cdot)|^2  \: + \:  \int_0^t \int_{F(s)}  2 \mu_F |D(u_F)|^2 ds \\
&  + \:  \frac{1}{2\beta_\Omega}\int_0^t \int_{\pa \Omega}  | u_F \times \nu|^2
\: + \:  \frac{1}{2\beta_S} \int_0^t  \int_{\pa S(t)}  |(u_F - u_S) \times \nu|^2    \\
& \le \: \int_0^t  \int_{F(t)} \rho_F (-g) \cdot u_F\: + \:  \int_0^t  \int_{S(t)} \rho_S (-g) \cdot u_S  \: + \:  \int_{\Omega \setminus S_0} \rho_F \, |u_{F,0}|^2  \: + \:  \int_{S_0} \rho_S \, |u_{S,0}|^2. 
\end{aligned}
\end{equation}
This goes together with the conservation of mass,  that amounts to the transport of $S$ by the rigid vector field $u_S$.  It reads 
$$ \pa_t \chi_S + \div (u_S \chi_S) \: = \:  0 \: \mbox{ in }  \Omega, \quad \chi_S(t,x) \: := \: 1_{S(t)}(x), $$
or in a weak form: for all $\Psi \in {\cal D}([0,T), {\cal D}(\overline{\Omega}))$,
\begin{equation} \label{mass}
 - \int_0^T \int_{S(t)} \pa_t \Psi \: - \: \int_0^T \int_{S(t)} u_S \cdot \na \Psi \: = \: \int_{S_0} \Psi\vert_{t=0}.
 \end{equation}

\noindent Pondering on these formal manipulations, we can now introduce our definition of a weak solution on $[0,T)$.  {\em We fix once for all the positive constants $\rho_S, \rho_F, \mu_F, \beta_S, \beta_\Omega$}. 
\begin{definition}
Let $\Omega$ and $S_0 \subset \Omega$ two Lipschitz bounded domains of $\R^3$. Let $u_{F,0} \in L^2_\sigma(\Omega)$, $u_{S,0} \in {\cal R}$
such that 
$\displaystyle
u_{F,0} \cdot \nu = u_{S,0} \cdot \nu \: \text{ on $\partial S_0$}.
$

\smallskip
\noindent
A weak solution of \eqref{NS}-\eqref{S3D}-\eqref{slip} on $[0,T)$  (associated to the initial data $S_0$, $u_{F,0}$ and  $u_{S,0}$) is a couple $(S, u)$ satisfying
\begin{itemize}
\item $S(t) \subset \Omega$  is a  bounded domain of $\R^3$ for all $t \in [0,T)$, such that 
$$\chi_S(t,x) := 1_{S(t)}(x)  \in L^\infty((0,T) \times \Omega). $$ 
\item   $u$ belongs to the space 
\begin{multline*}
{\cal S}_T \: := \: \Bigl\{ u \in L^\infty(0,T; L^2_\sigma(\Omega)), \: \mbox{there exists }  \: u_F \in L^2(0,T;  H^1_\sigma(\Omega)), \:  u_S \in L^\infty(0,T;  {\cal R})  \\
 \mbox{ such that } \:  u(t,\cdot) = u_F(t,\cdot) \:  \mbox{ on } F(t), \:  u(t,\cdot) = u_S(t,\cdot) \:  \mbox{ on } S(t), \: \mbox{ for a. e. } \:  t \in [0,T] \Bigr\},
\end{multline*}
where $F(t) := \Omega \setminus \overline{S(t)}$ for all $t \in [0, T)$. 
\item Equation \eqref{momentum} is satisfied  for all $\varphi \in {\cal T}_T$. 
\item Equation \eqref{mass} is satisfied for all $\psi \in {\cal D}([0,T); {\cal D}(\overline{\Omega}))$.
\item Equation \eqref{mainenergybound} is satisfied for almost every $t \in (0,T).$ 
\end{itemize}
\end{definition}
Let us conclude this paragraph by a few comments on this definition of weak solutions: 
\begin{enumerate}
\item As $\chi_S \in L^\infty((0,T) \times \Omega)$, the  integrals over $S(t)$ in \eqref{mass} are integrable with respect to time: namely, 
$$  t \mapsto \int_{S(t)} \pa_t \Psi = \int_\Omega \chi_S \pa_t \Psi  \quad   \mbox{ and }  \: t \mapsto \int_{S(t)} u_S \cdot \na \Psi = \int_\Omega \chi_S u_S \cdot \na \Psi $$
belong to $L^1(0,T)$. Actually, by the method of characteristics, as $u_S \in L^\infty(0,T; {\cal R})$ (rigid velocity field), it is easily seen that 
$$ S(t) = \phi_{t,0}(S_0) $$
for an isometric propagator $\phi_{t,s}$ which is Lipschitz continuous in  time, smooth in space. It follows that all integrals in equation \eqref{momentum} make sense. For instance, as  $\pa S(t)$ is Lipschitz for all $t$  and fields $u_F, u_S, \varphi_F, \varphi_S$ have  at least $L^2H^1$ regularity, the surface  integral 
$$ \int_{\pa S(t)}  ((u_F - u_S) \times n) \cdot ((\varphi_F - \varphi_S) \times n)  $$
 can be defined for almost every $t$ in the trace sense. Moreover, it defines an element  of  $L^1(0,T)$. This can be seen through the change of variable $x = \phi_{t,0}(y)$: the surface integral turns into 
$$    \int_{\pa S_0} \j\left(t,\phi_{t,0}(y)\right) \, \mbox{Jac}_\tau(y) \, dy, $$
  where 
  $$ \j(t,x) \: := \:      \left((u_F(t,x) - u_S(t,x)) \times \nu \right) \cdot \left((\varphi_F(t,x) -  \varphi_S(x)) \times \nu \right) $$
and where 
$$  \mbox{Jac}_\tau(y) \: = \: \| [\na \phi_{t,0}(y)]^{-1} \, \nu(y) \|_2  \det(\na \phi_{t,0}(y)) \: (=1) $$
is the tangential jacobian (see \cite[Lemme 5.4.1]{henrot} for details). This clearly defines an element of $L^1(0,T)$. 
\item Equations \eqref{momentum} and \eqref{mass} involve fields $u_F, u_S$, $\varphi_F,\varphi_S$ defined over $\Omega$ and such that 
$$ u = (1-\chi_S) u_F + \chi_S u_S, \quad \varphi = (1-\chi_S) \varphi_F + \chi_S \varphi_S,  $$
However, a closer look at equations \eqref{momentum} and \eqref{mass} shows that they only involve 
$$ \chi_S u_S,  \: \chi_F (1, \na) u_F,  \quad \mbox{ as well as } \: \chi_S (1, \pa_t) \varphi_S \:\mbox{  and } \:   \chi_F (1,\pa_t, \na) \varphi_F. $$
 In particular, they only depend on $u$ and $\varphi$, not on the choice of the extended fields $u_F, u_S$ and $\varphi_F, \varphi_S$. 
 \item The condition $u \in L^\infty(0,T; L^2_\sigma(\Omega))$ implies that 
 $$ u_F \cdot \nu\vert_{\pa S(t)} = u_S \cdot \nu\vert_{\pa S(t)} \: \mbox{ for a.e.} \: t   $$
 all terms being again understood in the trace sense. 
\item It is easy to see that equation \eqref{mass}, that is the  transport equation 
$$ \pa_t \chi_S + \div(\chi_S u_S)  = 0 \: \mbox{ in } \:  {\cal D}'([0,T) \times \overline{\Omega}) $$
can be written 
\begin{equation} \label{transportsolid}
 \pa_t \chi_S + \div(\chi_S u) = 0  \: \mbox{ in } \:  {\cal D}'([0,T) \times \overline{\Omega}) 
\end{equation}
and implies 
\begin{equation} \label{transportfluid}
 \pa_t \chi_F +  \div(\chi_F u) = 0  \: \mbox{ in } \:  {\cal D}'([0,T) \times \overline{\Omega}), \quad \chi_F(t,x) = \chi_{F(t)}(x) 
 \end{equation}
(remind that $F(t) = \Omega\setminus\overline{S(t)}$).  More generally, one can replace $u$ by any $v \in L^\infty(0,T; L^2_\sigma(\Omega))$ satisfying 
$$ v(t,\cdot) \cdot \nu\vert_{\pa S(t)} = u \cdot \nu \vert_{\pa S(t)} = u_S \cdot \nu \vert_{\pa S(t)} \quad \mbox{  for a.e.}\: t  $$
where the last equality holds in the space $H^{-1/2}(\pa S(t))$ (see \cite[Theorem 3.2.2]{Galdibooknew}). 
Note that equations \eqref{transportsolid} and \eqref{transportfluid} should be replaced by 
$$ \pa_t \rho_s + \div(\rho_s u) = 0, \quad  \pa_t \rho_f + \div(\rho_f u) = 0 $$
in the case of inhomogeneous solid and fluid,  with variable density functions $\rho_s$ and $\rho_f$. See \cite{Feireisl03} in the case of no-slip conditions. Extension of the present work (on a single  rigid and homogeneous solid in a homogeneous fluid) to more general configurations will be the matter of a forthcoming paper.  
\item Noticing that 
$$ D(u(t,\cdot)) = D(u_S(t,\cdot))  = 0 \: \mbox{ in } \: S(t), \quad D(\varphi(t,\cdot)) = D(\varphi_S(t,\cdot))  = 0 \: \mbox{ in } \: S(t) $$
it is tempting to write \eqref{momentum} under the condensed form 
\begin{multline} \label{momentum_unified}
 -\int_0^T \int_\Omega \rho u \cdot \partial_t\varphi + \int_0^T \int_\Omega \rho u \otimes u : D(\varphi) + \int_0^T \int_\Omega 2\mu_F D(u) : D(\varphi) \\
 \: = \: \mbox{ "boundary terms"} 
\end{multline}
where $\rho := \rho_F \chi_F + \rho_S \chi_S$, coupled with the global transport equation 
\begin{equation} \label{transport_unified}
\pa_t \rho + \div(\rho u) = 0 \: \mbox{ in } \: \Omega. 
\end{equation}
This kind of global formulation, reminiscent of the inhomogeneous Navier-Stokes equations, is used in the construction of weak solutions with Dirichlet boundary conditions: {\it cf} \cite{SanMartin&al02}. {\em However, it is not valid here}: due to the discontinuity of the tangential components of $u$ and $\varphi$,  neither $\pa_t \varphi$ nor $D(u)$ and $D(\varphi)$  belong to $L^2(\Omega)$. For instance, 
$$ \pa_t \varphi = \chi_F \pa_t \varphi_F \: + \:  \chi_S \pa_t \varphi_S \: + \: u_S \cdot \nu \left( \varphi_F - \varphi_S \right) \delta_{\pa S} $$   
where $\delta_{\pa S}$ is the Dirac mass along the solid boundary $\pa S$. This is why we keep the formulation \eqref{momentum}, distinguishing between the solid and the fluid part. 
\item  The definition of  a weak solution that we consider can not be satisfactory after collision. Indeed, we do not specify any rebound law. Moreover, in the case of Dirichlet conditions at the fluid-solid interface, explicit examples show that the analogue of our weak solution is not unique \cite{Starovoitov03bis}.  
\end{enumerate}

\subsection{Main result} \label{mainresult}
Our result is the following 
\begin{theorem} {\bf (Existence of weak solutions up to collision)} \label{theorem_existence}

\smallskip \noindent
Let $\Omega$ and $S_0 \Subset \Omega$ two $C^{1,1}$ bounded domains of $\R^3$. Let $u_{F,0} \in L^2_\sigma(\Omega)$, $u_{S,0} \in {\cal R}$
such that 
$\displaystyle 
u_{F,0} \cdot \nu = u_{S,0} \cdot \nu \: \text{ on $\partial S_0$}.
$

\smallskip
\noindent
There exists $T > 0$ and a  weak solution of \eqref{NS}-\eqref{S3D}-\eqref{slip} on $[0,T)$  (associated to the initial data $S_0$, $u_{F,0}$ and  $u_{S,0}$). Moreover,  such weak solution exists up to collision, that is 
$$ S(t) \Subset \Omega \quad \mbox{ for all } \:  t \in [0,T), \quad \mbox{ and } \: \lim_{t \rightarrow T^-} \mbox{dist}(S(t), \pa \Omega) = 0. $$   
\end{theorem}

The rest of the paper will be devoted to the proof of the theorem. Briefly, there are two main difficulties compared to the case of Dirichlet conditions:
\begin{itemize}
\item The lack of a unified formulation such as \eqref{momentum_unified}.
\item The lack of a uniform $H^1$ bound on  solutions $u$. 
\end{itemize}
These difficulties appear both in the construction of  approximate solutions, and in the convergence process. 

\medskip
Indeed, the approximation of fluid-solid systems is usually adressed by relaxing the solid constraint, through a penalization term. In this way, one is left with approximate systems that are close to density dependent Navier-Stokes equations. Roughly, they read 
\begin{equation} \label{approxpena}
\left\{
\begin{aligned}
&\pa_t (\rho_n u_n)  \: +  \:   \div(\rho_n u_n \otimes u_n)  \: +  \: \dots \: = \:  \mbox{ penalization} \\
& \pa_t \rho_n + \div(\rho_n u_n) = 0.
\end{aligned}
\right.
\end{equation}
In the case of no slip conditions, in which a global formulation of type \eqref{momentum_unified}-\eqref{transport_unified} already holds, to build such approximation is quite natural. But in the case of Navier conditions, this is not easy. 

\medskip
Once an approximate sequence of solutions $(\rho_n,u_n)$ has been derived,  Dirichlet conditions  allow for uniform $H^1$ bounds on $u_n$. This simplifies a lot of convergence arguments, notably with regards to the transport equation 
$$ \pa_t \rho_n + \div(\rho_n u_n) = 0 $$
to which the classical DiPerna-Lions theory applies straightforwardly \cite{DiPernaLions89}. Also, it helps to provide strong convergence of $u_n$ in $L^2((0,T) \times \Omega)$. In short, the lack of bound on $\pa_t (\rho_n u_n)$ (due to the penalization term) can be overcome by considering the fields $P_\delta(t) u_n$, where $P_\delta(t)$ is the orthogonal projection in $H^1_\sigma(\Omega)$ over the fields that are rigid in a $\delta$-neighborhood of $S(t)$. One can show that  $P_\delta(t) u_n$ has good equicontinuity properties uniformly in $\delta$ and $n$. 

\medskip
In the case of Navier boundary conditions, no uniform bound is available in $H^1$. This forces us to use more the structure of the solution $u$, in particular the $H^1$ bounds on the fluid and solid domains separately. This is also a source of trouble for the construction of approximate solutions, as one must find an approximation scheme in which such structure is not too much broken. 

\subsection{Strategy of proof} \label{strategy}
Let us describe here briefly the main lines of our proof. Let $S_0, u_{F,0}, u_{S,0}$ as in Theorem \ref{theorem_existence}, and 
$$ \rho_0 \: :=  \: \rho_F (1 - 1_{S_0})  \: + \: \rho_S 1_{S_0}, \quad u_0 \: :=  \:   (1 - 1_{S_0}) u_{F,0}  \: + \:  1_{S_0} u_{S,0}\,.$$
The keypoint is to consider approximate problems of the following type: {\em find $(S^n, u^n)$ such that 
\begin{description}
\item[a)] $\: S^n(t) \subset \Omega$ is a bounded Lipschitz domain for all $t \in [0,T]$, such that 
$$ \chi^n_S(t,x) \: := \: 1_{S^n(t)}(x) \in L^\infty((0,T) \times \Omega) \cap C([0,T];  L^p(\Omega)), \: \forall \, p < +\infty$$
\item[b)] 
$\displaystyle \:   u^n \in L^\infty(0,T; L^2_\sigma(\Omega)) \cap L^2(0,T; H^1_\sigma(\Omega)). $
\item[c)] For all $\varphi \in H^1(0,T;L^2_{\sigma}(\Omega)) \cap L^2(0,T; H^1_\sigma(\Omega))$ s.t. $\varphi\vert_{t=T} = 0$: 
\begin{align*}
& -\int_0^T \int_\Omega \rho^n \, \left(u^n \pa_t \varphi + v^n \otimes u^n : \na \varphi \right) \: + \:  \int_0^T \int_\Omega 2 \mu^n D(u^n) : D(\varphi) \\
&  + \:  \frac{1}{2\beta_\Omega}\int_0^T \int_{\pa \Omega}  (u^n \times \nu) \cdot (\varphi \times \nu)
\: + \:  \frac{1}{2\beta_S} \int_0^T  \int_{\pa S^n(t)}  ((u^n - P^n_S u^n) \times \nu) \cdot ((\varphi - P^n_S \varphi) \times \nu) \\ 
& \: + \:  n \int_0^T \int_\Omega  \chi^n_S  (u^n - P^n_S u^n)  \cdot (\varphi - P^n_S \varphi)   \: = \:  \int_0^T  \int_{\Omega} \rho^n (-g) \cdot \varphi  \: + \:  \int_{\Omega} \rho_0\, u_0 \cdot \varphi\vert_{t=0}
\end{align*} 
\item[d)] 
$ \displaystyle \: \pa_t \chi^n_S \: + \:  P^n_S u^n \cdot  \na \chi^n_S = 0, \quad \chi^n_S\vert_{t=0} = 1_{S_0}. $
\end{description}
}

\medskip
In above lines, 
\begin{itemize}
\item $\rho^n \: := \:   \rho_F (1 - \chi^n_S)  \: + \: \rho_S \chi^n_S$ is the total density function. 
\item $\mu^n \: := \: \mu_F (1 - \chi^n_S) \: + \: \frac{1}{n^2} \chi^n_S \:$ is an inhomogeneous  viscosity coefficient. 
\item $P^n_S = P^n_S(t)$  is the orthogonal projection in $L^2(S^n(t))$ over rigid fields. This means that: 
\end{itemize}
$$\forall \, 0 \le t < T, \quad \forall  u_S \in {\cal R}, \quad \forall u \in L^2_\sigma(\Omega), \quad P^n_S(t) u \in {\cal R} \quad \mbox{and } \quad \int_\Omega \chi^n_S(t, \cdot) (u - P^n_S(t) u) \, \cdot \,  u_S = 0. $$
\begin{itemize}
\item Eventually, 
$$v^n   \in L^\infty(0,T; L^2_\sigma(\Omega)) \cap L^2(0,T; H^1_\sigma(\Omega))$$
is a field that satisfies 
\begin{align*}
&v^n(t, \cdot) \: = \:  P^n_S(t) \, u^n(t, \cdot) \: \mbox{ in } \:  S^n(t), \\
& v^n(t, \cdot) \: = \:  u^n(t, \cdot) \: \mbox{ outside a   $ \delta$ neighborhood of } \: S^n(t), \: t \in [0,T),  
\end{align*}
for some $\delta$ fixed  and arbitrary in $\displaystyle (0, \mbox{dist}(S_0,\pa \Omega)/2)$. Moreover, $v^n$ will be chosen so that it is close to $u^n$ outside $S^n$ (in $L^p$ topology). In this way, it will asymptotically coincide with the limit $u$ of $u^n$. Further details on the definition of $v^n$ will be provided in due course.  
\end{itemize}

\medskip
Let us make a few comments on such approximate problems: 
\begin{enumerate}
 \item They rely on the use of the fields $P^n_S u^n$, that were already introduced in \cite{BostCottetMaitre10} in the context of Dirichlet conditions. These fields appear both: 
 
\smallskip
 i) in the transport equation for $\chi^n_S$. They will allow for a good control of the trajectories of the approximate solid bodies $S^n$. 
 
 \smallskip
 ii) in the penalization term    $\displaystyle n \int_0^T \int_\Omega  \chi^n_S  (u^n - P^n_S u^n)  \cdot (\varphi - P^n_S \varphi)$. Formally, as $n$ goes to infinity, this term will allow to recover the rigid constraint inside the solid. 
\item  Note that {\it a contrario} to the large penalization term, the viscosity term $\mu^n$ vanishes asymptotically in the solid part. Hence, there will be no uniform bound in $H^1_\sigma(\Omega)$ for $u^n$, as expected (see the discussion in  paragraph \ref{mainresult}). 
\item A specificity of these approximate problems is that the transport equation {\bf d)} is nonlinear in $\chi^n_S$ for a given $u^n$. Indeed, $P^n_S$ depends  on $\chi^n_S$ ({\it cf} the formula in section \ref{transport}). The whole section \ref{transport} is dedicated to this auxiliary nonlinear transport equation, which is a keystone of the approximation procedure.   
\item Once the solution $\chi^n_S$ of {\bf d)} is found and seen as a functional of $u^n$, equation {\bf c)} can be written as $\displaystyle {\cal F}(u^n) = 0$ for some functional ${\cal F}$  from $L^\infty L^2_\sigma \cap L^2 H^1_\sigma$ into itself. In short, we shall solve this equation by a Galerkin procedure: we shall look for an approximate solution  $u^{n,N}(t,x) = \sum_{k=0}^N \alpha_k(t) e_k(x)$ where $(e_k)$ is an orthonormal basis of $L^2_\sigma(\Omega)$. We shall solve approximate equations 
$\displaystyle {\cal F}^{n,N}(u^{n,N}) = 0$ by Schauder's theorem and pass to the limit with respect to $N$. This process is explained in section \ref{approx}.  
\item Note that the field $v^n$ satisfies
$$ v^n(t,\cdot)  \cdot \nu\vert_{\pa S^n(t)} \: = \: P^n_S(t) u^n(t,\cdot) \cdot \nu \vert_{\pa S^n(t)}. $$
 In particular, one can write 
$$ \pa_t \chi^n_S \: + \: v^n \cdot \na \chi^n_S = 0, \quad \mbox{ and } \quad \pa_t \rho^n    \: + \: v^n \cdot \na \rho^n = 0. $$
This will allow to obtain energy estimates in a standard way, in the spirit of the approximate systems \eqref{approxpena} used for Dirichlet conditions. The price to pay is the necessary  control of $u^n-v^n$, which will exhibit strong gradients near $\pa S^n$ as $n \rightarrow +\infty$. Moreover, a similar "boundary layer behaviour" will be involved in the approximation of discontinuous test functions $\varphi \in {\cal T}_T$ by continuous test functions $\varphi^n$ (involved in {\bf c)}). The whole convergence process will be analyzed in section \ref{convergence}. 
\end{enumerate}

 \section{A nonlinear transport equation} \label{transport}
Let $T > 0$, $u \in L^\infty(0,T; L^2_\sigma(\Omega))$. This section is devoted to the equation 
$$ \pa_t \chi_S \: + \: P_S u \cdot \na \chi_S, \quad \chi_S\vert_{t=0} = 1_{S_0}, $$ 
where $P_S u$ is defined by the following formula
\begin{equation} \label{defPs}
 P_S u \: :=   \: \frac{1}{M} \int_{\Omega} \rho_S \, \chi_S \, u \, \: + \: \left( J^{-1} \int_{\Omega} \rho_S  \, \chi_S  \, \left( (x'- x_S) \times u \right) dx') \right) \times (x- x_S)  
 \end{equation}
where the center of mass, total mass and inertia tensor of the solid are defined by 
\begin{equation} \label{eq_xS} x_S \: :=  \: \int_{\R^3} \rho_S \chi_S, \quad M \: := \: \int_{\R^3} \rho_S \chi_S, \end{equation}
and 
\begin{equation} \label{eq_J} J \: := \: \int_{\R^3} \rho_S \chi_S\left(  |x-x_S|^2 I_d -  (x-x_S) \otimes   (x-x_S) \right) dx.\end{equation}
If $\chi_S(t,x) = 1_{S(t)}(x)$ with $S(t)$ a subdomain of $\Omega$ , $P_S(t)$ is the orthogonal projection in  $L^2(S(t))$ over rigid vector fields, see \cite{BostCottetMaitre10}.  

\bigskip
We start with the regular case, that is  when  $ u \in C([0,T]; {\cal D}_\sigma(\overline{\Omega}))$. This case will be useful for Galerkin approximations of {\bf a)-d)}. 

\begin{proposition} {\bf (Well-posedness)} \label{transportWP}

\smallskip \noindent
Let $u \in  {C}([0,T]; {\cal D}_\sigma(\overline{\Omega})).$ 

\medskip \noindent
i) There is a unique solution $\:  \chi_S \in L^\infty((0,T) \times \R^3) \cap C([0,T]; L^p(\R^3))$ $(p < \infty) \: $ of 
\begin{equation} \label{transportR3}
\pa_t \chi_S + \div(\chi_S \, P_S u  ) = 0 \:  \mbox{ in } \R^3, \quad \chi_S\vert_{t=0} \: = \: 1_{S_0}. 
\end{equation}
\\
ii) Moreover $ \chi_S(t,\cdot) = 1_{S(t)}$ for all $t$, with   $S(t)$ a  Lipschitz bounded  domain. More precisely,   
$$ S(t) =  \phi_{t,0}(S_0)$$
 for the isometric propagator $\phi_{t,s}$ associated to $P_S u$ :   
$\displaystyle \: (t,s) \mapsto \phi_{t,s} \in C^1([0,T]^2; C^\infty_{loc}(\R^3))$.
\end{proposition}
{\em Proof.}
We can suppose  $u  \in  {C}([0,T]; {\cal D}_\sigma(\mathbb R^3))$ with no loss of generality. 

\medskip
Assume for a moment that we have found a solution $\chi_S$ of \eqref{transportR3}. Then, we can see \eqref{transportR3} as a linear transport equation, with given transport $P_S u \in C([0,T]; \: {\cal R})$. By the method of characteristics, we get easily
\begin{equation} \label{charac}   
 \chi_S(t,\phi_{t,0}(y)) = 1_{S_0}(y),  
 \end{equation} 
where $\phi_{t,s}$ is the isometric propagator defined by 
\begin{equation} \label{eq_phi}
\left\{
\begin{aligned}
\phi_{s,s}(y) &= y \,, \quad \forall \, y \in \mathbb R^3,\\
\partial_t \phi_{t,s}(y) &= P_S u(t,\phi_{t,s}(y) ) \quad \forall \, (s,t,x) \in (0,T)^2  \times \mathbb R^3. 
\end{aligned}
\right.
\end{equation}
Now, we  use \eqref{charac} in the expression \eqref{defPs} for $P_S u$. We obtain: 
\begin{multline} \label{eqPsu}
P_S u(t,x) = \frac{1}{M} \int_{S_0} \rho_S \, 1_\Omega(\phi_{t,0}(y)) \,  u(t,\phi_{t,0}(y)) dy  \\
+ \left( J^{-1}(t) \int_{S_0}  \rho_S \,  1_\Omega(\phi_{t,0}(y)) \,  \left(\phi_{t,0}(y) - x_S(t)\right) \times u(t,\phi_{t,0}(y)) dy \right) \times (x - x_S(t))
\end{multline}
where $\: M \: :=  \: |S_0| \, \rho_S, \quad  x_S(t) \: := \: \displaystyle{\int_{S_0}} \rho_S \phi_{t,0}(y)dy, \: $ and 
$$ J(t) \: := \: \int_{S_0} \rho_S \biggl(  | \phi_{t,0}(y)-x_S(t)|^2 I_d -  ( \phi_{t,0}(y)-x_S(t)) \otimes   ( \phi_{t,0}(y)-x_S(t)) \biggr) dy. $$
In particular, denoting ${\rm Isom}(\R^3) \approx \R^3 \times O_3(\R)$ the finite dimensional manifold of  affine isometries, we deduce from \eqref{eq_phi} and \eqref{eqPsu} that $t \mapsto \phi_{t,0}, \: [0,T] \mapsto {\rm Isom}(\R^3)$ satisfies an ordinary differential equation, of the type
\begin{equation} \label{odephi}
\frac{d}{dt} \phi_{t,0} = U_S(t, \phi_{t,0}), \quad \phi_{0,0} = I_d, 
\end{equation}
for a time-dependent vector field $U_S$ over ${\rm Isom}(\R^3)$. Namely, $U_S(t,\phi) \in T_\phi({\rm Isom}(\R^3)) \approx {\cal R}$ is defined by the same formula as in  \eqref{eqPsu}, replacing everywhere $\phi_{t,0}$ by $\phi$. 

\medskip
Conversely, if we manage to show existence of and uniqueness of a $C^1$ solution of \eqref{odephi} over $[0,T]$, then formula \eqref{charac} will define the unique solution $\chi_S$ of the nonlinear equation \eqref{transportR3}, proving Theorem \ref{transportWP}. 

\medskip
Hence, it only remains to study the well-posedness of \eqref{odephi}. We can identify ${\rm Isom}(\R^3)$ with $\R^3 \times O_3(\R) \subset \R^3 \times \R^9$, and identify all tangent spaces with ${\cal R}  \subset \R^3 \times \R^9$.   By the Cauchy-Lipschitz theorem, there is existence and uniqueness of a  $C^1$ maximal solution  if $U_S$ is continuous in $t,\phi$, locally Lipschitz in $\phi$.   Considering the expression of $U_S$, see \eqref{eqPsu}, this follows from 
\begin{lemma} \label{lem_M*}
Let $v \in C([0,T] ; {C}^{\infty}_{loc}(\R^3))$.  Then, the function 
$${\cal M}: [0,T] \times {\rm Isom}(\R^3) \mapsto \R, \quad   {\cal M}(t,\phi) = \int_{S^0} 1_{\Omega}(\phi(y)) v(t, \phi(y)) dy $$
 is continuous in $(t,\phi)$, and uniformly Lipschitz in $\phi$  over  $[0,T]$. 
\end{lemma}  
{\em Proof of the lemma.} The continuity is obvious. Then, for two affine isometries $\phi$ and $\phi'$, we write
\begin{align*}
 {\cal M}(t,\phi) - {\cal M}(t,\phi') & = \int_{S_0} 1_\Omega(\phi(y)) \, \left( v(t,\phi) - v(t,\phi'(y)) \right)  \: + \:   \int_{S_0}  \left( 1_{\Omega}(\phi(y)) -  1_{\Omega}(\phi'(y)) \right) \, v(t,\phi'(y)) \, dy \\
 & := \: M_1(t) + M_2(t). 
 \end{align*}
Clearly, 
$$ | M_1(t) | \: \le \:     \sup_{\substack{t \in [0,T],  \\ |x| \le  \|(\phi,\phi')\|_{\infty}}} \hspace{-.5cm}| \pa_x v(t,x)| \: \int_{S_0} | \phi(y) - \phi'(y)| dy  \: \le \: C_{\phi,\phi'} || \phi - \phi'||_{\infty}.  $$ 
 As regards $M_2$, we write 
 $$  M_2(t) \: \le \:   \sup_{\substack{t \in [0,T],  \\ |x| \leq \|\phi'\|_{\infty} }} \hspace{-.2cm}  |v(t,x)|   \:  \int_{\R^3}   \left| 1_{\Omega}(\phi(y)) -  1_{\Omega}(\phi'(y)) \right| \, dy \: \le \: C_{\phi'} \,  \int_{\R^3}  \left|  1_{\Omega}(\phi(y)) -  1_{\Omega}(\phi'(y)) \right| \, dy $$
For each $y$, the integrand is non-zero if and only if $\phi(y) \in \Omega$ and $\phi'(y) \in \Omega^c$ or vice-versa. As $|\phi(y) - \phi'(y) | \le  || \phi - \phi'||_{\infty}$, this is only possible if $\phi(y)$ and $\phi'(y)$ are in a $|| \phi - \phi'||_{\infty}$-neighborhood (say $V$) of $\pa \Omega$. Hence,
$$ | M_2(t)  | \: \le \: C_{\phi,\phi'} \, \left(  \int_{\phi^{-1}(V)}  dy +   \int_{\phi'^{-1}(V)}  dy \right) \: \le \:  2 C_{\phi,\phi'} | V | \: \le \:  C'_{\phi,\phi'}  || \phi - \phi'||_{\infty}.$$
This concludes the proof of the lemma.

\bigskip
Last step is to prove that the maximal solution is defined over the whole interval $[0,T]$.  From \eqref{eq_phi}-\eqref{eqPsu}, one can write
$$ \phi_{t,0}(y) = x_S(t) + Q_S(t) y $$
 where $x_S(t)$ is defined in \eqref{eq_phi} and $Q_S(t)$ is an orthogonal matrix. In particular, the only way that the maximal solution is not global on $[0,T]$ is through a blow-up of $x_S$. But, again, from \eqref{eqPsu},   
 $$ |\frac{d}{dt} x_S(t)| \: =  \: \frac{1}{M}  \left|   \int_{S_0} \rho_S \, 1_\Omega(\phi_{t,0}(y)) \,  u(t,\phi_{t,0}(y)) dy \right| \: \le  \: C \, \| u \|_{L^\infty((0,T)\times \Omega)}  $$
which prevents any blow-up. This ends the proof of the theorem. 
 
%
%

\begin{proposition}  {\bf (Strong sequential continuity)} \label{strongcontinuity}
\smallskip \noindent
Assume that   
$$\displaystyle  u^n \rightarrow  u \quad \mbox{  in } \: C([0,T];  {\cal D}_\sigma(\overline{\Omega})). $$
Then with obvious notations, one has 
\begin{equation*}
\chi^n_S \rightarrow \chi_S   \quad   \mbox{ weakly *  in } \: L^{\infty}((0,T) \times \R^3), \quad \mbox{ strongly in } \:   C([0,T];  L^p_{\alert loc}(\R^3)) \: ( p < \infty), 
\end{equation*}
as well as 
\begin{equation*}
 P^n_S u^n  \rightarrow P_S u \quad   \mbox{ strongly in } \:  C([0,T]; C^\infty_{loc}(\R^3)), \quad 
 \phi^n \rightarrow \phi \quad  \mbox{ strongly in } \: C^1([0,T]^2; C^\infty_{loc}(\R^3)).
\end{equation*} 
 \end{proposition}
 {\em Proof of the proposition.}
 As $u^n$ converges in $C([0,T];{\cal D}_{\sigma}(\overline{\Omega})),$  we have that $P^n_S u^n$
 is bounded in $L^{\infty}(0,T;L^2(\Omega))\,.$ Similarly $\chi^n_S$ is bounded in $L^{\infty}((0,T)\times \Omega).$
 Furthermore, up to  a subsequence that  we do not relabel, $P^n_S u^n$ converges weakly-* in $L^{\infty}(0,T;H^1_{loc}(\mathbb R^3))$ to some $\bar{u}_S$
 and $\chi^n_S(0) (= 1_{S_0}, \text{ for all $n \in \mathbb N$})$ converges strongly in $L^1(\Omega).$ Applying Di Perna-Lions theory, we obtain that 
 $\chi^n_S$ converges weakly-* in $L^\infty((0,T)\times \Omega)$ and strongly in $C([0,T];L^p_{loc}(\Omega))$ for all finite $p$. Its limit $\bar{\chi}_S$,  satisfies :
 $$
 \partial_t \bar{\chi}_S + \text{div} (\bar{\chi}_S  \bar{u}_S ) = 0\,. 
 $$ 
 Using  the convergence of both $\chi^n_S$ and  $u^n$ in equation \eqref{defPs}, we obtain that $\bar{u}_S = \bar{P}_S u$, where $\bar{P}_S$ is defined similarly to $P_S u$, replacing $\chi_S$ by $\bar{\chi}_S$. Moreover, the convergence of $P^n_S u$
 holds in $C([0,T]; C^\infty_{loc}(\R^3)).$ Consequently, $(\bar{\chi}_S, \bar{P}_S u)$ 
 is the unique solution of \eqref{transportR3} so that $\bar{\chi}_S = \chi_S$ and $\bar{P}_S u = P_S u$, and all the sequence converges.
 
 \medskip
To derive the convergence of the propagators $\phi^n$ from the  convergence of  the vector fields $P^n_S u^n$  is then standard, and we omit it for brevity.

 \begin{proposition}  {\bf (Weak sequential continuity)} \label{weakcontinuity}

\smallskip
\noindent
Let $(u^n, \chi^n_S)$ be a bounded sequence  in  $L^\infty(0,T; L^2_\sigma(\Omega)) \times L^\infty((0,T) \times \Omega)$, satisfying 
$$ \pa_t \chi^n_S + \div(P_S^n u^n \, \chi^n_S) = 0 \:  \mbox{ in } \R^3, \quad \chi^n_S\vert_{t=0} \: = \: 1_{S_0}. $$
Then, up to a subsequence,  one has 
$$ u^n \rightarrow u \quad \mbox{ weakly *  in } \:   L^\infty(0,T; L^2_\sigma(\Omega)) $$
\begin{equation*}
 \chi^n_S \rightarrow \chi_S \quad  \mbox{ weakly *  in } \: L^{\infty}((0,T) \times \R^3), \quad  \mbox{ strongly in} \:   C([0,T];  L^p_{loc}(\R^3)) \: ( p <  \infty), 
 \end{equation*}
with  $(u_S, \chi_S)$ a solution of 
$$ \pa_t \chi_S + \div(P_S u \, \chi_S) = 0 \:  \mbox{ in } \R^3, \quad \chi_S\vert_{t=0} \: = \: 1_{S_0}.  $$
Moreover,  $\chi_S$ satisfies condition ii) of  Proposition \ref{transportWP}, and the following additional convergences hold: 
\begin{eqnarray*}
 P^n_S u^n  \rightarrow P_S u \quad &   \mbox{ weakly * in } &  \:  L^\infty(0,T; C^\infty_{loc}(\R^3)),  \\[6pt]
\phi^n \rightarrow \phi  \quad  &\mbox{ weakly *  in } & \: W^{1,\infty}((0,T)^2; C^\infty_{loc}(\R^3))  \quad  \mbox{ strongly in} \:   C([0,T]^2;  C^\infty_{loc}(\R^3))\,.
\end{eqnarray*}
 \end{proposition}
 {\em Proof.} 
 The proof follows the same scheme as the previous one. We only sketch the arguments.
 First,  up to the extraction of a subsequence, we obtain that 
 $$
  u^n \rightarrow u \quad \mbox{ weakly *  in } \:   L^\infty(0,T; L^2_\sigma(\Omega))
 $$
Then, as before, we obtain  that $P^n_S u^{n}$ is bounded in  $L^{\infty}(0,T; C^\infty_{loc}(\R^3))$. This yields that 
$$
P_S^n u^n \rightarrow \bar{P}u  \mbox{ weakly *  in } \:   L^\infty(0,T;  H^1_{loc}(\mathbb R^3)) 
$$
(still up to a subsequence). We then deduce applying Di Perna-Lions theory that, up to the extraction of a subsequence,  $\chi^n_S$ converges strongly in $ C([0,T];  L^p_{loc}(\R^3))$
to some $\chi_S$,  which in turn implies that $\bar{u}_S = P_S u$ and that $(\chi_S, P_S u)$ is a solution to our tranport equation. Eventually,  uniform bounds on $\phi^n$ and $\pa_t \phi^n$ (which imply weak-* convergence of a subsequence in $W^{1,\infty}$) follow easily.  

\section{Approximation} \label{approx}
This section is devoted to the resolution of approximate fluid-solid systems. These approximate systems were introduced in paragraph \ref{strategy}, {\it cf} {\bf a)-d)}. The previous section has focused on the transport equation {\bf d)}. It remains to examine {\bf c)}. At first, we  explain a little how the field  
$v^n$ connecting $P^n_S u^n$ to ${u}^n_S$ is defined.   The detailed definition of $v^n$ will be achieved in section \ref{convergence}. 

\subsection{Connecting velocity} \label{mollified}
We first remind a classical result on the equation $\div u = f$, taken from \cite[Exercise III.3.5]{Galdibooknew}:
\begin{proposition} \label{Bogovski}
Let ${\cal O}$  be a bounded Lipschitz domain. Let  $f \in L^2({\cal O})$ and 
$\varphi \in H^{1/2}(\pa{\cal O})$ satisfying the compatibility condition
$$ \int_{{\cal O}} f = \int_{\pa {\cal O}} \varphi \cdot \nu. $$
 Then there exists a solution $u \in H^1({\cal O})$ of
$$ \div u = f \quad \mbox{ in } \:  {\cal O}, \quad u = \varphi \quad \mbox{ at } \:  \pa {\cal O} $$ 
with 
$$ \| u \|_{H^1(\cal O)} \: \le \: C_{{\cal O}} \left( \| \varphi \|_{H^{1/2}(\partial \cal O)} \: + \: \| f \|_{L^2(\cal O)} \right).  $$ 
\end{proposition}
The previous proposition yields easily  

\begin{corollary} \label{lemmamollified} {\bf (Extension of solenoidal vector fields)}

\smallskip
\noindent
There exists a continuous linear operator $E_\Omega : H^1_{\sigma}(\overline{\Omega}) 
\mapsto H^1_\sigma(\R^3)$
satisfying $E_\Omega \, u  = u$ on $\Omega$. Moreover, for all open subset $\omega \Subset \Omega$, 
$$ \| E_\Omega \, u \|_{H^1(\R^3 \setminus \overline{\omega})} \: \le \: C_\omega  \| u \|_{H^1(\Omega \setminus \overline{\omega})}, \quad \forall \: u \in H^1_\sigma(\overline{\Omega}). $$
\end{corollary}

\begin{corollary} {\bf (Connection  of solenoidal vector fields)} \label{coroconnection}

\smallskip
\noindent
For all $\delta > 0$, there exists a continuous linear operator 
 $$ V^\delta :  H^1_\sigma\left(\R^3\setminus S_0\right) \times   H^1_\sigma(\overline{S_0}) \, \mapsto \, H^1_\sigma(\R^3), \quad  (U, U_S) \, \mapsto \,  V^\delta[U, U_S]  
 $$
 such that  
 \begin{align*}
&V^\delta[U,U_S] \: = \:  U_S \quad \mbox{ in   } \:  S_0, \\
& V^\delta[U,U_S] \: = \: U \quad  \mbox{ outside a   $\delta$ neighborhood of } \: S_0.   
\end{align*}
\end{corollary}
From there, we have the following 
\begin{proposition} \label{constructionvn}
 For all $\delta > 0$, there exists a continuous mapping 
 $$ v^\delta : L^2(0,T; H^1_{\sigma}(\R^3)) \times   L^\infty(0,T; {\cal R}) \, \mapsto \,  L^2(0,T; H^1_{\sigma}(\R^3)), \quad (u, u_S) \mapsto  v^\delta[u, u_S]  $$
such that 
\begin{align*}
&v^\delta[u,u_S](t, \cdot) \: = \:  u_S(t, \cdot) \quad \mbox{ in  } \:  S(t), \\
& v^\delta[u,u_S](t, \cdot) \: = \:  u(t, \cdot) \quad \mbox{ outside a   $\delta$ neighborhood of } \: S(t), \: t \in [0,T),  
\end{align*}
where, as usual, $\: S(t) := \phi_{t,0}(S_0)$ and  $\: \phi = \phi_{t,s}$ is the isometric propagator associated to $u_S$. Moreover, $v^\delta$ can be chosen so that  
 $$ \| v^\delta[u,u_S] \|^2_{L^2(0,T; H^1(\mathbb R^3))} \: \le \: {\cal C} \, \int_0^T \left( \| u(t,\cdot) \|^2_{H^1(\R^3\setminus \overline{S(t)})} \: + \:   \| u_S(t,\cdot) \|^2_{L^2(S(t))} \right) \, dt, $$
 where ${\cal C}$ depends on $\delta$ and  $T$.  
\end{proposition}

\medskip
{\em Proof of the proposition.} The proposition can be deduced from Corollary \ref{coroconnection} using Lagrangian coordinates. Namely, we introduce $U$ and $U_S$ through the relations
$$ u\left(t,\phi_{t,0}(y)\right) \: = \: d\phi_{t,0}\vert_{y}\left(U(t,y) \right) , \quad u_S\left(t,\phi_{t,0}(y)\right) \: = \: d\phi_{t,0}\vert_{y}\left(U_S(t,y) \right).    $$
Clearly, for all $t$, $\: U(t,\cdot)$ and $U_S(t,\cdot)$ define elements of $H^1_{\sigma}(\R^3\setminus S_0)$ and $H^1_\sigma(\overline{S_0})$ respectively. 
Using Corollary \ref{coroconnection}, we  define $v^\delta[u,u_S]$ through the relation  
$$  v^\delta[u,u_S]\left(t,\phi_{t,0}(y)\right) \: = \: d\phi_{t,0}\vert_{y}\left(V^\delta[U(t,\cdot),U_S(t,\cdot)](y) \right).  $$ 
It fulfills all requirements, which ends the proof.

 \bigskip
 Back to system {\bf c)}, the idea is to define 
 $$ v^n \: :=  \:   \: v^{\delta}[E_\Omega u^n, P^n_S u^n] . $$
Clearly, for any time $T^n$ such that 
 $$ \mbox{dist}(S^n(t), \pa \Omega) \ge 2 \delta, \quad t \in [0,T^n], $$
 $v^n\vert_{\Omega}$ will belong to $L^2(0,T^n; H^1_\sigma(\Omega))$ and will satisfy  
 \begin{align*}
&v^n(t, \cdot) \: = \:  P^n_S(t) \, u^n(t, \cdot) \quad \mbox{ in  } \:  S^n(t), \\
& v^n(t, \cdot) \: = \:  u^n(t, \cdot) \quad \mbox{ outside a   $\delta$ neighborhood of } \: S^n(t), \: t \in [0,T^n).  
\end{align*}
Let us stress that there is still some latitude left in the construction of $v^n$, through the  choice of the operator $V^\delta$ in Corollary \ref{coroconnection}. As will be shown  in section \ref{convergence}, this operator can be chosen depending on $n$ ($V^\delta = V^{\delta,n}$) so that 
$v^n$ is close to  $u^n$  outside $S^n$ (in $L^p$ topology). However, this additional property will not be needed until section \ref{convergence}.

\medskip
Last remark:  the resolution of {\bf a)-d)}, and the whole construction of weak solutions, will be first performed on a small time interval $[0,T]$, for a time $T$ that is uniform in $n$. Existence of weak solutions up to collision will follow from a continuation argument, to be explained at the end of section \ref{convergence}.

\subsection{Galerkin approximation}
As pointed out in paragraph \ref{strategy}, the resolution of  {\bf a)-d)} is carried out through a Galerkin scheme. Let $(e_k)_{k \ge 1}$ being both  an orthonormal basis of $L^2_\sigma(\Omega)$ and a basis of $H^1_\sigma(\Omega)$,  with elements in ${\cal D}_\sigma(\overline{\Omega})$. The aim of this paragraph is to find for all $N, n$ and some $T > 0$  a couple $(S^N, u^N)$ satisfying
\begin{description}
\item[a')] $\: S^N(t) \subset \Omega$ is a bounded Lipschitz domain for all $t \in [0,T]$, such that 
$$ \chi^N_S(t,x) \: := \: 1_{S^N(t)}(x) \in L^\infty((0,T) \times \Omega) \cap C([0,T];  L^p(\Omega)) \:  (p < \infty) $$
\item[b')] 
$\displaystyle \:   u^N(t,\cdot) \: = \: \sum_{i=1}^N \alpha_k(t) e_k, \quad \mbox{ with } \: \alpha = (\alpha_1, \dots, \alpha_N) \in C([0,T])^N.$
\item[c')] For all $\varphi \in {\cal D}([0,T); \: \mbox{ span}(e_1, \dots, e_N))$
\begin{align*}
& -\int_0^T \int_\Omega \rho^N \,  \left(u^N \cdot \pa_t \varphi + v^N \otimes u^N : \na \varphi \right) \: + \:  \int_0^T \int_\Omega 2 \mu^N D(u^N) : D(\varphi) \\
&  + \:  \frac{1}{2\beta_\Omega}\int_0^T \int_{\pa \Omega}  (u^N \times \nu) \cdot (\varphi \times \nu)
\: + \:  \frac{1}{2\beta_S} \int_0^T  \int_{\pa S^N(t)}  ((u^N - P^N_S u^N) \times \nu) \cdot ((\varphi - P^N_S \varphi) \times \nu) \\ 
& + n \int_0^T \int_\Omega  \chi^N_S  (u^N - P^N_S u^N)  \cdot (\varphi - P^N_S \varphi)   \: = \:  \int_0^T  \int_{\Omega} \rho^N (-g) \cdot \varphi  \: + \:  \int_{\Omega} \rho_0\, u_0 \cdot \varphi\vert_{t=0}
\end{align*} 
\item[d')] 
$ \displaystyle \: \pa_t \chi^N_S \: + \:  P^N_S u^N \cdot  \na \chi^N_S = 0 \: \mbox{ in } \: \Omega, \quad \chi^N_S\vert_{t=0} = 1_{S_0}. $
\end{description}
In above lines, similarly to the original problem:
\begin{itemize}
\item $\rho^N \: := \:   \rho_F (1 - \chi^N_S)  \: + \: \rho_S \chi^N_S$ is the total density function. 
\item $\mu^N \: := \: \mu_F (1 - \chi^N_S) \: + \: \frac{1}{n^2} \chi^N_S \:$ is an inhomogeneous  viscosity coefficient. 
\item $P^N_S = P^N_S(t)$  is defined by \eqref{defPs}, adding the upperscript $N$ everywhere. 
\item Eventually, $v^N = v^{\delta}[u^N, P^N_S u^N]$, see paragraph \ref{mollified}.   
\end{itemize}
{\em Note that all quantities above depend on $n$, notably through the penalization term and the viscosity coefficient. But we omit $n$ from the notations to lighten writings.} Also, note that $u^N$ can be seen as an element of $L^2(0,T; H^1_\sigma(\R^3))$, as the $e_k$ are defined globally.  In particular,  $v^N = v^{\delta}[u^N, P^N_S u^N]$ is well-defined. 

\medskip
The main result of this paragraph is 
\begin{theorem} \label{theogalerkin}
There is $T > 0$, $R > 0$, such that for all $n,N$, {\bf a')-d')} has  at least one solution such that $\| u^N \|_{L^\infty(0,T; L^2(\Omega))} \: \le \:  R$. 
\end{theorem}
To prove Theorem \ref{theogalerkin}, we shall express our Galerkin problem as a fixed point problem, and will apply Schauder's theorem to it. Thus, we want to identify $u^N$ as the fixed point of an application 
$$ {\cal F}^N : u \mapsto \tilde u, $$
defined on $B_{R,T} \: := \: \left\{ u  \in C([0,T]; \: \mbox{span}(e_1,\dots, e_N)), \quad  \| u \|_{L^\infty(0,T; L^2(\Omega))} \le R \right\}$. 
We proceed as follows. Let $u\in B_{R,T}$. 
\begin{itemize}
\item {\em Step 1.} Let $\chi_S$ 
be the solution of 
$$ \displaystyle \: \pa_t \chi_S \: + \:  P_S u \cdot  \na \chi_S = 0, \quad \chi_S\vert_{t=0} = 1_{S^0}, $$
given by Proposition \ref{transportWP}. We know that $\chi_S(t,x) = 1_{S(t)}(x)$ with $S(t)$ a bounded Lipschitz domain, $t \in [0,T]$. We  define accordingly: 
$$ \rho \: := \:   \rho_F (1 - \chi_S)  \: + \: \rho_S \chi_S, \quad \mu \: := \:  \mu_F (1 - \chi_S) \: + \: \frac{1}{n^2} \chi_S, \quad v(t,x) \: :=  \: v^\delta[u, P_S u]. $$
\item {\em Step 2.} We consider the following ODE, with unknown $\tilde u : [0,T] \mapsto \mbox{ span}(e_1,\dots, e_N)$: 
\begin{equation} \label{odegalerkin}
 A(t) \frac{d}{dt} \tilde u(t)  + B(t) \tilde u(t) = f(t), \quad \tilde u(0) \: = \:  u^N_0 \: := \: \sum_{k=1}^N \left( \int_{\Omega} u_0 \cdot e_k \right) \, e_k, 
 \end{equation}
in which  $A(t) :=  (a_{i,j}(t))_{1 \le i,j \le N}$, $B(t) := (b_{i,j}(t))_{1 \le i,j \le N}$  and $f(t) := (f_i(t))_{1\le i \le N}$ are defined by 
\begin{align*}
a_{i,j} \: & := \:  \int_{\Omega} \rho e_i \cdot e_j , \\
b_{i,j} \: & := \:   \int_\Omega \rho (v \cdot \na e_j) \cdot e_i  \: + \:   \int_\Omega 2 \mu D(e_i) : D(e_j) \: + \: 
 \frac{1}{2\beta_\Omega}\int_{\pa \Omega}  (e_i \times \nu) \cdot (e_j \times \nu) \\ 
& + \:  \frac{1}{2\beta_S}  \int_{\pa S(t)}  ((e_i - P_S e_i) \times \nu) \cdot ((e_j - P_S e_j) \times \nu) \: + \:  n \int_{\Omega}  \chi_S  (e_i - P_S e_i)  \cdot (e_j - P_S e_j)  \\
 f_i \: & := \: \int_{\Omega} \rho (-g) \cdot e_i.
\end{align*}
We have identified here the function $\tilde u$ with its coefficients in the basis $e_1, \dots, e_N$.  Note that the function $\rho$ defined in step 1  has a positive lower bound, so that $ A(t) \geq \min(\rho_S,\rho_F)I_N$ in the sense of symmetric matrices, whatever the value of  $\chi_S$. Also, the continuity of $A$ and $B$ over $[0,T]$ is easy and will be proved below. In particular, equation \eqref{odegalerkin} has a unique solution 
$$ \tilde u \in C^1([0,T];  \: \mbox{span}(e_1,\dots, e_N)). $$
\end{itemize}
In this way, we can associate to each $u \in B_{R,T}$ some field  
$$ \tilde u = {\cal F}^N(u) \in  C([0,T];  \: \mbox{span}(e_1,\dots, e_N)). $$
The whole point is to prove 
\begin{proposition} \label{propgalerkin}
There exists $T > 0$, $R > 0$, uniform in $n$ and $N$, such that 
${\cal F}^N$ is a  well-defined mapping from $B_{R,T}$ to itself, continuous and compact. 
\end{proposition}
Before proving this proposition, let us show how it implies Theorem \ref{theogalerkin}. By Schauder's theorem, it yields the existence of a fixed point $u^N \in B_{R,T}$ of ${\cal F}^N$. Let  $\chi^N_S = 1_{S^N}$ be the corresponding solution of the transport equation on $[0,T] \times \R^3$. As will be clear from the proof, the time  $T$  of the proposition satisfies     
$$ \mbox{dist}(S^N(t), \pa \Omega) \: \ge \: 2 \delta, \quad \forall \, t \in [0,T], $$
for some $\delta$ fixed  and arbitrary in $\displaystyle ( 0, \mbox{dist}(S_0,\pa \Omega)/2)$. 
Hence,  {\bf a')} is satisfied,  and $v^N := v^{\delta}[u^N, P^N_S u^N]$ satisfies  $\: v^N \cdot \nu \vert_{\pa \Omega} = 0$, as well as  
$$ \pa_t \rho^N  + v^N \cdot \na \rho^N = 0 \: \mbox{ in } \Omega $$
(see  remark 4 after the definition of weak solutions, and remark 5, paragraph \ref{strategy}). Finally, we notice  that ODE \eqref{odegalerkin} is equivalent to:  for all $\varphi \in {\cal D}([0,T); \: \mbox{ span}(e_1, \dots, e_N))$

\begin{equation}
\label{unifestim}
\begin{aligned}
&  \int_0^T\int_\Omega \rho^N \,  \pa_t  u^N \cdot  \varphi \: + \: \int_0^T\int_\Omega \rho^N v^N \cdot \na  u^N  \cdot \: \varphi \: + \:    \int_0^T \int_\Omega 2 \mu^N D(u^N) : D(\varphi) \\
&  + \:   \frac{1}{2\beta_\Omega}  \int_0^T \int_{\pa \Omega}  (u^N \times \nu) \cdot (\varphi \times \nu)
\: + \:  \frac{1}{2\beta_S}   \int_0^T   \int_{\pa S^N(t)}  ((u^N - P^N_S u^N) \times \nu) \cdot ((\varphi - P^N_S \varphi) \times \nu) \\ 
& + n  \int_0^T \int_\Omega  \chi^N_S  (u^N - P^N_S u^N)  \cdot (\varphi - P^N_S \varphi)   \: = \:   \int_0^T \int_{\Omega} \rho^N (-g) \cdot \varphi 
- \int_{\Omega}  \rho_0 \,  u_0^N \varphi\vert_{t=0} \,.
\end{aligned}
\end{equation}  
Combining this equation with the previous one on $\rho^N$ leads to {\bf c')}. Note that condition $\: v^N \cdot \nu \vert_{\pa \Omega} = 0$ is needed for  the convective term to vanish through integration by parts.

\medskip
{\em Proof of the proposition.}

\medskip
{\em Step 1:  Definition of ${\cal F}^N$}. 

\medskip
We first prove that ${\cal F}^N$ is well-defined from $B_{R,T}$ to $C([0,T];  \: \mbox{span}(e_1,\dots, e_N))$ for any $T$ and $R > 0$. The only thing to check is the continuity of matrices $A$ and $B$ in \eqref{odegalerkin} with respect to time, which will guarantee the existence of a solution to the linear ODE \eqref{odegalerkin}. As $\chi_S$ belongs to  $C([0,T]; L^p(\Omega))$ for all finite $p$, so does $\rho$, and  $A$ is clearly continuous.  As regards $B$, the only difficult terms are 
$$ I(t) \: := \:   \int_\Omega \rho (v(t,\cdot) \cdot \na e_j) \cdot e_i, \quad J(t) \: := \:  \frac{1}{2\beta_S}  \int_{\pa S(t)}  ((e_i - P_S(t) e_i) \times \nu) \cdot ((e_j - P_S(t) e_j) \times \nu). $$
We remind that the propagator $\phi = \phi_{t,s}$ associated to $P_S u$ satisfies 
$$ \phi \in C^{1}\left([0,T]^2; \, C^\infty_{loc}(\R^3)\right) $$
Hence, a look at the construction of $v^\delta$, {\it cf} Corollary \ref{coroconnection} and Proposition \ref{constructionvn} (see also Lemma \ref{lemmaw} in the appendix \ref{app_lemmawn}), yields 
$$ v   \in C([0,T]; H^1_\sigma(\R^3)).$$
 It implies that  $t \mapsto I(t)$ is continuous. 
 
 \medskip
 As regards $J(t)$, we change variables to go back to a fixed domain. We set $x = \phi_{t,0}(y)$ to obtain 
 $$ J(t) \: :=  \:  \frac{1}{2\beta_S}  \int_{\pa S_0} \j\left(t,\phi_{t,0}(y)\right) \, \mbox{Jac}_\tau(y) \, dy, $$
  where 
  $$ \j(t,x) \: := \:      \left((e_i(x) - P_S(t) e_i(x)) \times \nu\right) \cdot \left((e_j(x) - P_S(t) e_j(x)) \times \nu\right) $$
and where 
$$  \mbox{Jac}_\tau(y) \: = \: \| d\phi_{t,0}\vert^{-1}_y \, \nu(y) \|_2  \det(d\phi_{t,0}\vert_y) (=1) $$
is the tangential jacobian. See \cite[Lemme 5.4.1]{henrot} for details. 
As $\j$ is continuous in $t$ and smooth in $x$, we obtain that $t \mapsto J(t)$ is continuous. 

\medskip
{\em Step 2: ${\cal F}^N$ sends $B_{R,T}$ to itself}. 

\medskip
Here, we need to restrict to small $T$. More precisely, we fix $0 < \delta < \frac{1}{2} \mbox{dist}(S_0, \pa \Omega)$, and consider a time $T$ such that 
\begin{equation} \label{conditionnoncontact}
 \inf_{u \in B_{R,T}} \, \mbox{dist}(S(t), \pa \Omega) \ge 2 \delta > 0 
\end{equation}  
{\em Let us prove  that such time $T$ does exist and can be chosen uniformly with respect to $N$ and $n$. } For all $u \in B_{R,T}$, we  write
$$ S(t) \: = \: \phi_{t,0}(S_0) $$
with $\phi$  the propagator associated to the rigid field $\: P_S u =  \dot{x}_S + \omega_S \times (x-x_S)$ defined in \eqref{defPs}. 
It is enough that 
$$ \sup_{t \in [0,T]} | \pa_t \phi_{t,0}(t,y) |< \frac{\mbox{dist}(S_0, \pa \Omega) - 2 \delta}{T}, \quad t \in [0,T], \quad y \in S_0. $$
We find
$$  | \pa_t \phi_{t,0}(t,y) | \: < \: |u_S(t, \phi_{t,0}(t,y)) | <   |\dot{x}_S(t)| \: + \: |\omega_S(t)| \, |y - x_{S_0}| $$  
using that the propagator is isometric. Moreover, classical  calculations yield
$$ | \dot{x}_S(t) |^2 + J(t) \, \omega_S(t) \cdot \omega_S(t) = \int_{S(t)} \rho_S |P_S u(t,\cdot)|^2 \le \int_{S(t)} \rho_S |u(t,\cdot)|^2 \: \le \:  \rho_S \, R^2.$$
We  can then use  the  inequality 
\begin{align*} 
 |\dot{x}_S(t)| \: + \: |\omega_S(t)| \, |y - x_{S_0}|  \: & \le \:  \sqrt{2} \, \max(1, |y - x_{S_0}|)  \,   \left( | \dot{x}_S(t) |^2 + |\omega_S(t)|^2 \right)^{1/2} \\
  & \le \: C_0 \,  \left(  | \dot{x}_S(t) |^2 + J(t) \omega_S(t) \cdot \omega_S(t) \right)^{1/2}
  \end{align*}
 where for instance 
 $$ C_0 \: := \:  \sqrt{2} \, \frac{ \max\left (1, \sup_{y \in S_0} |y - x_{S_0}|\right)}{\min(1, \lambda_0)^{1/2}}, \quad \lambda_0  \: : \: \mbox{
 smallest eigenvalue of $J(0)$}. $$ 
Eventually, any $\: \displaystyle  T \: < \: \frac{\mbox{dist}(S_0, \pa \Omega) - 2 \delta}{C_0 (\rho_S)^{1/2} R}$   will satisfy \eqref{conditionnoncontact}. 

\medskip

Let now $u$ be arbitrary in $B_{R,T}$. 
Thanks to \eqref{conditionnoncontact},  we have that
$\displaystyle  v = v^{\delta}[u, P_S u] $
 satisfies $v  \cdot \nu\vert_{\pa \Omega} = 0$, and   
$$ \pa_t \rho  + v \cdot \na \rho = 0 \: \mbox{ in } \Omega. $$
Mutiplying \eqref{odegalerkin} by $\tilde{u}$, integrating in time, and combining with  the last transport equation, we obtain the  energy estimate
\begin{align} \label{eq_energyN}
& \| \sqrt{\rho} \tilde u(t, \cdot) \|^2_{L^2} \: + \:  \int_0^t \int_\Omega  2 \mu |D(\tilde u)|^2 \notag \\
&  + \:  \frac{1}{2\beta_\Omega}\int_0^t \int_{\pa \Omega}  |\tilde u \times \nu|^2
\: + \:  \frac{1}{2\beta_S} \int_0^t  \int_{\pa S(t)}  |(\tilde u - P_S \tilde u) \times \nu|^2    + n \int_0^t \int_\Omega  \chi_S  |\tilde u - P_S  \tilde u|^2  \\
& \le \: \int_0^t  \int_{\Omega} \rho (-g) \cdot \tilde u  \: + \:  \int_{\Omega} \rho_0\, | u^N_0|^2 \notag 
\end{align}
As $\min(\rho_F,\rho_S) \le \rho \le   \max(\rho_F,\rho_S)$, we deduce easily that 
$$ \| \tilde u \|_{L^\infty(0,T; L^2(\Omega))} \: \le \: R $$
for $R = R(T, u_0)$ large enough. Hence, ${\cal F}$ sends $B_{R,T}$ to itself.

\medskip
{\em Step 3. Compactness of ${\cal F}^N$.}

\medskip
For any $\displaystyle u = \sum_{k=1}^N \alpha_k e_k$,  we get from  equation \eqref{odegalerkin}:
$$ | \frac{d}{dt} \tilde \alpha(t) | \: \le \:  | A^{-1}(t)| \, |B(t) | \, | \alpha(t) | \: + \: | f(t)| \: \le \:  \, R \,  | A^{-1}(t)| \, |B(t) |  \: + \: | f(t)|. $$
Integrating with respect to time, we obtain 
$$ \sup_{t \in [0,T]} \left( |\tilde\alpha(t)| + |\frac{d}{dt}\tilde\alpha(t)| \right)   \: \le \: C' $$
(where the constant at the r.h.s. may depend on $N$ or $n$). In other words, 
$$ \sup_{u \in B_{R,T}} \| {\cal F}^N(u) \|_{C^1([0,T]; \: \mbox{span}(e_1,\dots, e_N))} \: \le \: C'' $$
which provides compactness in  $B_{R,T}$ by Ascoli's theorem.

\medskip
{\em Step 4. Continuity of ${\cal F}^N$.}  

\medskip
Let $(u^k)$ a sequence in $B_{R,T}$, such that $u^k \rightarrow u$ in  $B_{R,T}$ (that is uniformly over $[0,T]$). We want to show that $\displaystyle {\cal F}^N(u^k) \rightarrow    {\cal F}^N(u)$ in $B_{R,T}$.  First, we note that, as $span(e_1,\ldots,e_N)$ is a finite-dimensional subspace of $\mathcal{D}_{\sigma}(\bar{\Omega})$ we have that $u^k$ converges to $u$ in $C([0,T] ; \mathcal{D}_{\sigma}(\bar{\Omega})).$ Then, we use Proposition \ref{strongcontinuity}. With obvious notations, 
 \begin{equation*}
\chi^k_S \rightarrow \chi_S   \quad   \mbox{ weakly *  in } \: L^{\infty}((0,T) \times \R^3), \quad \mbox{ strongly in } \:   C([0,T];  L^p_{loc}(\R^3)) \: ( p < \infty), 
\end{equation*}
as well as 
\begin{equation*}
 P^k_S u^k  \rightarrow P_S u \quad   \mbox{ strongly in } \:  L^\infty(0,T; C^\infty_{loc}(\R^3)), \quad 
 \phi^k \rightarrow \phi \quad  \mbox{ strongly in } \: W^{1,\infty}((0,T)^2; C^\infty_{loc}(\R^3)).
\end{equation*} 
From there, and the construction of $v^\delta$ (Corollary \ref{coroconnection}, Proposition \ref{constructionvn}, Lemma \ref{lemmawn} in appendix \ref{app_lemmawn}), it is easy to see that
$$ v^k \rightarrow v \: \mbox{ strongly in }  \: C([0,T]; H^1_{\sigma}(\R^3)). $$
 By slightly adapting the arguments of Step 1, one can then show that the matrices in \eqref{odegalerkin} satisfy 
 $$ B^k \rightarrow  B, \quad A^k \rightarrow A   \quad  \mbox{ strongly in }  \: C([0,T]).  $$
 From classical results for ODE's, it follows that
 $$ \tilde u^k = {\cal F}(u^k) \rightarrow  \tilde u = {\cal F}(u)   \: \mbox{ strongly in }  \: C([0,T]; \: \mbox{span}(e_1,\dots, e_N)). $$
 For the sake of brevity, we leave the details to the reader.

 \subsection{Convergence of the Galerkin scheme}
 In the previous paragraph, we have built for each $n, N$ a solution $u^{n,N}$ (denoted $u^N$ for brevity) of {\bf a')-d')}. It is defined on $[0,T]$ for some  time $T$ uniform in $n,N$, satisfying \eqref{conditionnoncontact}. The next step is to let $N$ go to infinity, to recover a solution $u^n$ of {\bf a)-d)}. We remind the uniform energy estimate  (see \eqref{unifestim})
  \begin{align} & \| \sqrt{\rho^N} u^N(t, \cdot) \|^2_{L^2(\Omega)} \: + \:  \int_0^t \int_\Omega  2 \mu^N |D(u^N)|^2 \: + \:  \frac{1}{2\beta_\Omega}\int_0^t \int_{\pa \Omega}  | u^N \times \nu|^2 \notag \\
&\: + \:  \frac{1}{2\beta_S} \int_0^t  \int_{\pa S^N(t)}  |(u^N - P^N_S u^N) \times \nu|^2    + n \int_0^t \int_\Omega  \chi^N_S  |u^N - P_S^N  u^N|^2  \label{estimationN}\\ 
& \le \: \int_0^t  \int_{\Omega} \rho (-g) \cdot u^N  \: + \:  \int_{\Omega} \rho_0\, |u^N_0|^2 \notag
\end{align}
It yields that 
$$ (u^N)_{N \in \N} \quad \mbox{ is bounded uniformly with respect to $N$ in } \: L^\infty(0,T; L^2_\sigma(\Omega)) \cap  L^2(0,T; H^1_\sigma(\Omega)) $$
The bound in $H^1$ follows from the $L^2$ bound on $D(u^N)$ and Korn's inequality, see \cite{Nitsche:1981}. From there, we will be able to show strong convergence both in the transport equation {\bf d')} and in the momentum equation {\bf c')}. As regards the transport equation, we rely on Proposition 
\ref{weakcontinuity}. Up to a subsequence, one has 
$$ u^N \rightarrow u \quad \mbox{weakly * in } \:   L^\infty(0,T; L^2_\sigma(\Omega)) \: \mbox{and weakly in} \:  L^2(0,T; H^1_\sigma(\Omega)) $$
for some $u (= u^n)$, and it follows from this proposition that
\begin{equation*}
 \chi^N_S \rightarrow \chi_S \quad  \mbox{ weakly *  in } \: L^{\infty}((0,T) \times \R^3), \quad  \mbox{ strongly in} \:   C([0,T];  L^p_{loc}(\R^3)) \: ( p <  \infty), 
 \end{equation*}
 as well as 
\begin{eqnarray*}
 P^N_S u^N  \rightarrow P_S u \quad &   \mbox{ weakly * in } &  \:  L^\infty(0,T; C^\infty_{loc}(\R^3)),  \\[6pt]
\phi^N \rightarrow \phi  \quad  &\mbox{ weakly *  in } & \: W^{1,\infty}((0,T)^2; C^\infty_{loc}(\R^3))\,,  \quad  \mbox{ strongly in} \:   C([0,T];  C^\infty_{loc}(\R^3))\,.
\end{eqnarray*}
%
up to another extraction. We stress again that all limits depend on $n$. 

\medskip
It remains to study the convergence of equation  {\bf c')}. Therefore, we fix the test function: we take 
$$\varphi(t,x) :=  \chi(t) \, e_j, \quad \chi \in {\cal D}([0,T)). $$
 for some fixed $j$. The point is to obtain as $N \rightarrow +\infty$ the limit  equation {\bf c)}, still with $\varphi(t,x) = \chi(t) e_j(x)$. But as $j$ is arbitrary, and as $(e_k)_{k \ge 1}$ is a basis of $H^1_\sigma(\Omega)$, standard density arguments will allow to extend the formulation to general test functions.   

\medskip
At first, we need to prove that, 
$$ v^N\vert_{\Omega} \: \rightarrow \:  v = v^{\delta}[E_\Omega u, P_S u] \vert_{\Omega} \quad \mbox{ in } \: L^2(0,T;H^1_\sigma(\Omega)). $$
It is enough to prove that 
 $$\hat{v}^N \: := \:  v^{\delta}[E_\Omega u^N, P^N_S u^N]  \:  \rightarrow \: \hat{v} \: := \:  v^{\delta}[E_\Omega u, P_S u]  $$
 weakly in $L^2(0,T; \, H^1_{loc}(\R^3))$. In view of Corollary \ref{coroconnection} and Proposition \ref{constructionvn}, it is an easy consequence of 
 Lemma \ref{lemmawn}.
 
 \medskip
 We are now ready to handle the asymptotics of {\bf c')} (with $\varphi(t,x) =  \chi(t) \, e_j(x)$). As before, for the sake of brevity, we focus on the two most difficult terms, those which involve  
\begin{eqnarray*}
 I^N(t) \:& := & \:   \int_\Omega \rho^N (v^N \otimes  u^N) : \nabla e_j, \\
J^N(t) \: &  :=  & \:  \frac{1}{2\beta_S}  \int_{\pa S^N(t)}  ((u^N - P^N_S(t) u^N) \times \nu) \cdot ((e_j - P^N_S(t) e_j) \times \nu). 
\end{eqnarray*}

 As regards $J^N(t)$, once again we change variables to go back to a fixed domain. We obtain 
 $$ J^N(t) \: :=  \:  \frac{1}{2\beta_S}  \int_{\pa S_0} \j^N\left(t,\phi^N_{t,0}(y)\right) \, \mbox{Jac}^N_\tau(y) \, dy, $$
  where 
  $$ \j^N(t,x) \: := \:      \left((u^N(t,x) - P^N_S(t) u^N(t,x)) \times \nu\right) \cdot \left((e_j(x) - P^N_S(t) e_j(x)) \times \nu\right) $$
and where 
$$  \mbox{Jac}^N_\tau(y) \: = \: \| d\phi^N_{t,0}\vert^{-1}_y \, \nu(y) \|_2  \det(d\phi^N_{t,0}\vert_y) = 1.  $$
Let $r^N \: :=  \: E_\Omega u^N - P^N_S u^N$, resp. $\eta_j^N := e_j - P^N_S e_j$ to which we associate $R^N$, resp. $H_j^N$ 
through the change of coordinates: 
$$
r^N(t,\phi^N(t,y)) := d\phi^N_t\vert_y R^N(t,y) \,, \qquad 
\eta_j^N(t,\phi^N(t,y)) := d\phi^N_t\vert_y H_j^N(t,y) \,.
$$ 
From  the weak convergence of $u^N$, we deduce that $r^N$ converges weakly in $L^2(0,T;H^1_{loc}(\mathbb R^3)).$
Given the strong convergence of  $\chi^N$ in $C([0,T];L^p(\Omega))$ we also have that $\eta_j^N$ converges strongly to $\eta_j := e_j - P_S e_j$ 
in $L^2(0,T;H^1_{loc}(\mathbb R^3))$. Furthermore, as $d\phi^N_{t,0}\vert_{y}$ is an isometric mapping for all $N$, we get that :
$$
\j^N(t,\phi_{t,0}(y)) = (R^N \times \nu)  \cdot  (H_j^N \times \nu )\,, \quad \forall \, N \in \mathbb N\,
$$
where, because of lemma \ref{lemmawn} :
$$
R^N \to  R \: \text{ weakly in } \: L^2(0,T;H^1(\Omega))\,, \quad H^N_j \to H_j \: \text{ strongly in } \: L^2(0,T;H^1(\Omega))\,.
$$
with obvious notations. This yields corresponding weak and strong convergences of the traces of these functions on $\partial S_0.$ Having in mind that $\mbox{Jac}^N_\tau \equiv 1$ for all $N$, and going back to the moving domain,  
we obtain easily that $J^N$ converges weakly in $L^1(0,T)$ to :
$$
J(t) := \frac{1}{2\beta_S}  \int_{\pa S(t)}  ((u - P_S(t) u) \times \nu) \cdot ((e_j - P_S(t) e_j) \times \nu). 
$$ 
\medskip

We finally turn to the convergence of $I^N,$ for which we will  need some  compactness on $(\rho^N u^N).$ Therefore, we introduce some notations:
 we denote by $P$ the orthogonal projection from $L^2(\Omega)$ onto $L^2_{\sigma}(\Omega),$ respectively $P_k$ the orthogonal projection 
from $L^2(\Omega)$ onto $\text{span}(e_1,\ldots,e_k)$. We also remind that  our strong, resp. weak, convergence results on $\rho^N$, resp. $u^N$ imply that.   
$$ \rho^N u^N \rightarrow \rho u \: \mbox{  weakly-* in } \: L^{\infty}(0,T; L^2(\Omega)). $$
In particular, we have for any fixed $k$:  
\begin{equation} \label{convergence1}
P_k(\rho^N u^N) \rightharpoonup  P_k(\rho u ) \quad \mbox{  weakly-* in } \: L^{\infty}(0,T; L^2_\sigma(\Omega)) \:  \mbox{as} \:  N \to \infty.
 \end{equation}
 Moreover,  equation {\bf c')}  can be written: for all $1 \le k  \le N$,  
$$\partial_t \, P_k(\rho^N u^N) + P_k F^N  = 0 \quad  \mbox{in} \quad   {\cal D}'\left(0,T; [H^1_\sigma(\Omega)]^*\right)$$
 where  $F^N \in {\cal D}'\left(0,T;  [H^1_\sigma(\Omega)]^*)\right)$ 
   is defined by the duality relation:   
 \begin{eqnarray*}
\langle F^N , \varphi \rangle  &=& 
 \: \int_0^T\int_\Omega \rho^N v^N \otimes  u^N :  \na \varphi \: - \:    \int_0^T \int_\Omega 2 \mu^N D(u^N) : D(\varphi) \\
& + &  \frac{1}{2\beta_S}   \int_0^T   \int_{\pa S^N(t)}  ((u^N - P^N_S u^N) \times \nu) \cdot ((\varphi - P^N_S \varphi) \times \nu) \\ 
& + & n  \int_0^T \int_\Omega  \chi^N_S  (u^N - P^N_S u^N)  \cdot (\varphi - P^N_S \varphi)  + \:   \frac{1}{2\beta_\Omega}  \int_0^T \int_{\pa \Omega}  (u^N \times \nu) \cdot (\varphi \times \nu)\\
& - &   \int_0^T \int_{\Omega} \rho^N (-g) \cdot \varphi, \qquad \mbox{ for all $\varphi \in \mathcal{D}\left(0,T ; H^1_\sigma(\Omega)\right).$}
\end{eqnarray*}
We remind that for $\displaystyle f \in [H^1_\sigma(\Omega)]^*$, $P_k$ is defined by duality: $<P_k f, \varphi> \: := \: < f, P_k \varphi>$.   
From the above expression for $F^N$  and the various bounds already obtained, it is easily seen that for any fixed $k$, $(P_k F^N)$ is bounded (in $N$) in  $\displaystyle L^2(0,T; [H^1_{\sigma}(\Omega)]^*)$.  Hence, the same conclusion applies to  $(\partial_t \, P_k(\rho^N u^N))$. Combining  with \eqref{convergence1}, it follows that for any fixed $k$,
\begin{equation} \label{convergence2}
P_k(\rho^N u^N) \rightarrow  P_k(\rho u ) \quad \mbox{  strongly in } \: L^{\infty}(0,T; [H^1_{\sigma}(\Omega)]^*) \:  \mbox{as} \:  N \to \infty.
 \end{equation}
Now, we note that, for arbitrary $k$ and $N,$ and a.a. $t \in (0,T)$ there holds 
\begin{eqnarray*}
\|P(\rho^N u^N )(t) - P_k  (\rho^N u^N)(t)\|_{[H^1_{\sigma}(\Omega)]^* }  &=& \sup_{\|\varphi\|_{[H^1_{\sigma}(\Omega)]} = 1} \int_{\Omega} \left(P(\rho^N u^N(t) ) - P_k  (\rho^N u^N)(t)]\right)
 \varphi\\
								&=& \sup_{\|\varphi\|_{[H^1_{\sigma}(\Omega)] } = 1} \int_{\Omega} \rho^N u^N(t)  (\varphi - P_k \varphi) \\
								&\leq & \left[  \sup_{\|\varphi\|_{[H^1_{\sigma}(\Omega)]} = 1} \|\varphi - P_k \varphi\|_{L^2(\Omega)} \right] \|\rho^N u^N\|_{L^\infty L^2(\Omega)}\,,\\
\end{eqnarray*}
By a standard argument based on Rellich Lemma, one shows that  
$$ \sup_{\|\varphi\|_{H^1_{\sigma}(\Omega)} =1} \|\varphi - P_k \varphi\|_{L^2(\Omega)} \rightarrow 0
$$ 
as $k \to \infty$. With the uniform bound on $\rho^N u^N$ in $L^\infty(0,T; L^2(\Omega))$, we can conclude that    
\begin{equation} \label{convergence3} 
P_k(\rho^N u^N) - P(\rho^N u^N) \: \rightarrow 0 \:  \quad \mbox{ strongly in  } \: L^{\infty}(0,T; [H^1_{\sigma}(\Omega)]^*), \quad \mbox{ as $k \rightarrow +\infty$, uniformly in $N$.}
\end{equation}
Of course, with a similar but simpler estimate, we also have 
\begin{equation} \label{convergence4} 
P_k(\rho u) - P(\rho u) \: \rightarrow 0 \:  \quad \mbox{ strongly in  } \: L^{\infty}(0,T; [H^1_{\sigma}(\Omega)]^*), \quad \mbox{ as $k \rightarrow +\infty$.} 
\end{equation}
Combining \eqref{convergence2}, \eqref{convergence3} and \eqref{convergence4}, we obtain finally that :
$P(\rho^N u^N)$ converges to $P(\rho u)$ strongly in $L^2(0,T;  [H^1_{\sigma}(\Omega)]^*).$
Combining this strong convergence with the weak convergence of $(u^N)$ in $L^2(0,T ; H^1_{\sigma}(\Omega)),$
we might apply the method of P.L. Lions \cite[p.47]{PLLvol1} with the duality bracket $[H^1_{\sigma}(\Omega)]^*-H^1_{\sigma}(\Omega) $ to prove
that $\sqrt{\rho^N} u^N$ converges to $\sqrt{\rho}  u$ strongly in $L^2((0,T) \times \Omega)$.
Finally, we rewrite :
$$
I_N(t) = \int_{\Omega} \sqrt{\rho^N} u^N \otimes \sqrt{\rho^N} v^N : \nabla e_j\,,  
$$ 
where :
\begin{itemize}
\item $\sqrt{\rho^N} u^N$ converges to $\sqrt{\rho}  u$ strongly in $L^2((0,T) \times \Omega)$
\item $\sqrt{\rho^N}$ converges to $\sqrt{\rho}$ strongly in $L^{\infty}(0,T ; L^3(\Omega))$
\item $v^N$ converges to $v$ weakly in $L^2(0,T ; L^6(\Omega))$ (thanks to the imbedding $H^1(\Omega) \subset L^6(\Omega)$). 
\end{itemize}
Combining these statements, we get that $I_N$ converges to $I$ (with obvious notations) weakly in $L^1(0,T).$ 

Such convergences result yield that $(\rho^n,u^n)$ satisfy {\bf c')} for test functions $\varphi$ of the form $\chi(t) \psi$ with $\chi \in \mathcal{D}([0,T))$ and $\psi \in {\rm span}(\{e_k, k \in \N\}).$
Via a classical density argument, the convergence extends to all $\displaystyle \varphi \in H^1(0,T; L^2_{\sigma}(\Omega)) \cap L^2(0,T;  H^1_\sigma(\Omega))$
such that $\varphi\vert_{t=T} = 0$.

\medskip

%

\subsection{Energy inequality} \label{sec_nrjin}
We end this section by proving that the approximate solution $(\rho^n,u^n)$ satisfies the further estimate :
\begin{align} & \| \sqrt{\rho^n} u^n(t, \cdot) \|^2_{L^2(\Omega)} \: + \:  \int_0^t \int_\Omega  2 \mu^n |D(u^n)|^2 \: + \:  \frac{1}{2\beta_\Omega}\int_0^t \int_{\pa \Omega}  | u^n \times \nu|^2 \notag \\
&\: + \:  \frac{1}{2\beta_S} \int_0^t  \int_{\pa S^N(t)}  |(u^n - P^n_S u^n) \times \nu|^2    + n \int_0^t \int_\Omega  \chi^n_S  |u^n - P_S^n  u^n|^2 \label{energyboundn} \\ 
& \le \: \int_0^t  \int_{\Omega} \rho (-g) \cdot u^n  \: + \:  \int_{\Omega} \rho_0\, |u_0|^2 \notag
\end{align}
for almost all $t \in [0,T].$ For simplicity we drop exponent $n$ in what follows. 

First, we note that the solutions $(\rho^N,u^N)$ of the Galerkin scheme satisfy \eqref{estimationN} uniformly in $N$ and that, up to the extraction of a subsequence $\sqrt{\rho^N}u^N$ converges to $\sqrt{\rho} u$ in $L^2((0,T) \times \Omega).$ Hence, we may pass to the limit in \eqref{estimationN} for almost
all $t \in [0,T]\,.$ On the other hand, there holds:
\begin{itemize}
\item By construction of the Galerkin scheme, $u_0^N \to u_0$ in $L^2(\Omega)$ so that :
$$
\lim_{N \to \infty} \int_{\Omega} \rho_0 |u_0^N|^2 = \int_{\Omega} \rho_0 |u_0|^2 \,.
$$
\item Given the strong convergence of $(u^N)$ in $L^2((0,T) \times \Omega)$ :
$$
\lim_{N \to \infty} \int_0^T \int_{\Omega} \rho(-g) \cdot u^N  = \int_0^T \int_{\Omega} \rho(-g) \cdot u^n\,.
$$
\item Given the weak convergence of $u^N$ in $L^2((0,T);H^1(\Omega))$ and the strong convergence of $\chi^N_S$ in $C([0,T];L^p(\Omega))$
we get that $\sqrt{\mu^N} D(u^N)$ converges weakly to $\sqrt{\mu} D(u)$ in $L^{2-\varepsilon}((0,T) \times \Omega).$ In particular, in the weak limit, there holds :
$$
\int_{0}^T \int_{\Omega} \mu |D(u)|^2 \leq \liminf \int_0^T \int_{\Omega} \mu^N |D(u^N)|^2 \,.
$$
With similar arguments, we obtain :
$$
\int_0^T \int_{\Omega} \chi_S |u -P_Su|^2 \leq \liminf \int_0^T \int_{\Omega} \chi_S^N |u^N - P^N_S u^N|^2 \,.
$$
\item Finally, we pass to the limit in the boundary terms. First, we introduce $U^N$  and $U^N_S$  associated to the extension $E_{\Omega}[u^N]$ and 
the rigid vector field $P^N_S u^N$ respectively, computed through the change of variable $\phi^N_{t,0}.$ As previously, we have:
$$
\int_{0}^T \int_{\partial S^N(t)} |(u^N - P_S^N u^N) \times \nu|^2  = \int_{0}^T \int_{\partial S^N(t)} |(u^N -  u_S^N) \times \nu|^2 =
\int_0^T \int_{\partial S_0} |(U^N - U_S^N) \times \nu|^2.
$$ 
Because of the weak convergence of $u^N$ and $u_S^N$ in $L^2(0,T; H^1_{\sigma}(\Omega)),$ we have that  
$U^N$ and $U_S^N$ converge also weakly in $L^2(0,T;H^1_{loc}(\R^3))$ (see Lemma \ref{lemmawn}). Hence, we have also weak convergence of the traces
on $S_0.$ The lower semi-continuity of the $L^2$-norm on $\partial S_0$ yields :  
\begin{eqnarray*}
\int_{0}^T \int_{\partial S(t)} |(u - P_S u) \times \nu|^2  &=& \int_{0}^T \int_{\partial S_0} |(U - U_S) \times \nu|^2  \\
	& \leq & \liminf \int_{0}^T \int_{\partial S_0} |(U^N - U^N_S) \times \nu|^2\\
	& \leq & \liminf  \int_{0}^T \int_{\partial S^N(t)} |(u^N - P^N_Su^N) \times \nu|^2
\end{eqnarray*}
Similar weak-convergence and semi-continuity arguments yield also :
$$
\int_0^T \int_{\partial \Omega} |u \times \nu|^2 \leq \liminf \int_0^T \int_{\partial \Omega} |u^N \times \nu|^2\,.
$$
\end{itemize}
This ends the proof of \eqref{energyboundn}.

 \section{Convergence} \label{convergence}
In the previous section, we have obtained the existence of solutions $u^n$ of approximate fluid-solid systems, namely satisfying {\bf a)}-{\bf d)}. These solutions $u^n$ are defined on some (uniform in $n$) time interval $(0,T)$ such that 
\begin{equation} \label{deltafar}
 \mbox{dist}(S^n(t), \pa \Omega) \ge 2 \delta, \quad \mbox{ for } \: t \in [0,T), \quad  \mbox{ for some fixed } \: \delta > 0. 
\end{equation}
We must now study the asymptotics of $u^n$  as $n$ goes to infinity, and recover a weak solution at the limit. 

\medskip
In what follows, we will often make use of the notation 
$$ ({\cal O})_\eta \: := \: \{ x \in \R^3, \:  \mbox{dist}(x, {\cal O}) < \eta \} $$
for ${\cal O}$ an open set and $\eta > 0$. 

\subsection{A priori bounds on $u^n$. Convergence in the transport equation} \label{aprioribounds}
The density  $\rho^n$ clearly  satisfies the uniform bound
\begin{equation} \label{boundrhon}
 \min(\rho_F,\rho_S) \le \rho^n \le  \max(\rho_F,\rho_S).
 \end{equation}
 Combining \eqref{energyboundn} and \eqref{boundrhon} yields  that
 \begin{equation} \label{boundun}
  \| u^n \|^2_{L^\infty(0,T;L^2(\Omega))}  + n \|  \sqrt{\chi^n_S} (u^n - P_S^n u^n) \|^2_{L^2((0,T) \times \Omega)}  +  \|\sqrt{\mu^n} D(u^n) \|^2_{L^2((0,T) \times \Omega)} \le C.
  \end{equation}
 for some constant $C$ depending only on $\rho_F, \rho_S$, $u_0$ and $\: T$. 
 
 \medskip
 In particular, up to a subsequence, the first inequality gives
 $$ u^n \rightarrow u \quad \mbox{ weakly* in } \: L^\infty(0,T; L^2_\sigma(\Omega)). $$ 
We can then pass to  the limit of the transport equation {\bf d)}, using Proposition \ref{weakcontinuity}. 
The following convergence holds up to a subsequence: 
\begin{equation*}
 \chi^n_S \rightarrow \chi_S \quad  \mbox{ weakly *  in } \: L^{\infty}((0,T) \times \R^3), \quad  \mbox{ strongly in} \:   C([0,T];  L^p_{loc}(\R^3)) \: ( p <  \infty), 
 \end{equation*}
with
$$\chi_S(t,\cdot) =  1_{S(t)}, \quad   S(t) = \phi_{t,0}(S_0)$$
 for an isometric propagator $\displaystyle \phi = \phi_{t,s} \in W^{1,\infty}((0,T)^2; C^\infty_{loc}(\R^3))$. Moreover, one has 
\begin{equation*}
 P^n_S u^n  \rightarrow P_S u \quad   \mbox{ weakly *  in } \:  L^\infty(0,T; C^\infty_{loc}(\R^3)),  \quad 
\phi^n \rightarrow \phi  \quad  \mbox{ weakly *  in } \: W^{1,\infty}((0,T)^2; C^\infty_{loc}(\R^3)).
\end{equation*}
In particular, one recovers the transport equation \eqref{mass}  setting $u_S := P_S u$. 

\medskip
Now,  we can combine  the second inequality in \eqref{boundun}, that  yields 
$$ \chi^n_S (u^n - P_S^n u^n) \: \rightarrow 0 \: \mbox{ strongly in } \: L^2,$$
with the strong (resp. weak) convergence of $\chi^n_S$ (resp. $u^n$ and $P^n_S u^n$). As $n$ goes to infinity,  we derive easily: 
\begin{equation} \label{solidedansS}
\chi_S  (u -  u_S)  = 0. 
\end{equation}

\medskip
Finally,  the last inequality in \eqref{boundun} and Korn's inequality imply that 
$$  \int_0^T  \|  u^n \|^2_{H^1(F^n(t))} dt \: \le \:  C \int_0^T \left(    \|  D(u^n) \|^2_{L^2(F^n(t))}  \: + \:  \|  u^n \|^2_{L^2(F^n(t))} \right) dt \: \le \: C, \quad F^n(t) \: := \: \Omega\setminus \overline{S^n(t)}.  $$
We then introduce continuous extension operators 
$$E_n(t) : \: \{ u \in H^1(F^n(t)), \: \div u = 0 \: \mbox{ in } F^n(t) , \: u \cdot \nu\vert_{\pa \Omega} = 0 \}  \: \mapsto \: H^1_\sigma(\Omega), $$
in the spirit of Corollary \ref{lemmamollified}.  As long as the $S^n(t)$ are $2 \delta$ away from $\pa \Omega$, it is standard to construct these extension operators in such a way that 
$$ \| E_n(t) \|_{{\cal L}(H^1)} \le C_\delta, \quad \forall \: t \in [0,T]. $$
Hence, if we set $u^n_F(t,\cdot) \: := E_n(t) u^n(t,\cdot)$, we have that 
$$ (u^n_F) \quad  \mbox{ is bounded in } \: L^2(0,T; \,  H^1_\sigma(\Omega)), \quad (1 - \chi^n_S) (u^n_F - u^n) = 0,  \: \forall \: n. $$
From the $L^2(0,T; H^1(\Omega))$ bound, we can assume up to another extraction that 
$$ u^n_F \rightarrow u_F \quad \mbox{weakly in }  \: L^2(0,T; H^1_\sigma(\Omega)). $$
From above equality and from the strong convergence of $\chi^n_S$, we then get:
\begin{equation} \label{fluidedansF}
(1 - \chi_S) (u_F - u) \: = \:  0. 
\end{equation}

\medskip
Eventually, considering relations \eqref{solidedansS} and \eqref{fluidedansF}, we get that the limit $u$ of $u^n$ belongs to ${\cal S}_T$. Hence, back to the definition of a weak solution, it only  remains to show that  the momentum equation \eqref{momentum} is satisfied by $S(\cdot), u_S, u_F$.

\subsection{A priori bounds on $v^n$.}  \label{aprioriboundsv}
Prior to  the analysis of the momentum equation, we must establish some refined estimates on the {\em connecting velocity} $v^n$. We remind that $v^n$
was defined in Lagrangian like coordinates, see paragraph \ref{mollified}. More precisely, 
$$ v^n\left(t, \phi^n_{t,0}(y)\right) \: := \: d\phi^n_{t,0}\vert_{y}\left( V^\delta[U^n(t,\cdot), U^n_S(t,\cdot)] \right),  $$
where 
$$ E_\Omega u^n\left(t,\phi^n_{t,0}(y)\right) \: = \: d\phi^n_{t,0}\vert_{y}\left(U^n(t,y) \right) , \quad P^n_S u^n\left(t,\phi^n_{t,0}(y)\right) \: = \: d\phi^n_{t,0}\vert_{y}\left(U^n_S(t,y) \right),  $$
and 
$V^\delta = V^\delta[U,U_S]$ is some linear operator connecting $\displaystyle \:  U \in H^1_\sigma(\R^3\setminus S_0) \: $ to $ \: \displaystyle  U_S \in H^1_\sigma\bigl(\overline{S_0}\bigr)$ over a band of width $\delta$ outside $S_0$: see Corollary \ref{coroconnection}. We shall here specify our choice for the operator $V^\delta$. Actually, we shall make it depend on $n \:$ ($V^\delta = V^{\delta,n}$), in order for the following additional assumption to be satisfied: 
\begin{multline} \label{extra}
\| V^{\delta,n}[U,U_S] - U \|_{L^p((S_0)_\delta \setminus S_0)} \le C_{\delta,p} \left( \| (U - U_S) \cdot \nu \|_{L^p(\pa S_0)} +  n^{1/6-1/p}     \| (U, U_S)\|_{H^1((S_0)_\delta \setminus S_0)} \right),\\  \quad \forall \,2 \le p \le 6. 
\end{multline}
We postpone  the construction of such operator $V^{\delta,n}$ to the end of the paragraph. 

\medskip
Back to $v^n$, the additional assumption \eqref{extra} implies easily  that for all $2 \le p \le 6$, 
\begin{multline}
 \| (1- \chi^n_S) (v^n -  u^n) \|_{L^2(0,T; L^p(\Omega))} \\
  \le \: C_{\delta,p} \left( \int_0^T \| (u^n - P^n_S u^n)(t,\cdot) \cdot \nu \|_{L^p(\pa S^n(t))}^2 \, dt \: + \: n^{1/6-1/p} \right). 
  \end{multline}
But we know that 
\begin{align} \notag  
\int_0^T \| (u^n - P^n_S u^n)(t,\cdot) \cdot \nu \|_{H^{-1/2}(\pa S^n(t))}^2 \, dt  \:& \le \:  C \,  \int_0^T \| (u^n - P^n_S u^n)(t,\cdot) \cdot \nu \|^2_{L^2(S^n(t))} \, dt\\ & \le \frac{C}{n} \label{unn1}
\end{align}
where the last bound comes from the second inequality in \eqref{boundun}. We emphasize here that, as the $S^n(t)$'s are all isometric to one another, 
the constant $C$ does not depend on $n,t$.  Interpolation with the similar other bound
\begin{multline}
\int_0^T \| (u^n - P^n_S u^n)(t,\cdot) \cdot \nu \|_{H^{1/2}(\pa S^n(t))}^2 \\
\begin{array}{rrcl}
&\: \le \: &    C \displaystyle{\int_0^T} \left( \| u^n \|_{H^1( (S^n(t))_{\delta} \cap F^n(t))}^2  + \:  \| u^n \|^2_{L^2(S^n(t))} \right) \text{d$t$}\, \\[4pt]
&\:\le\: &  C \displaystyle{\int_0^T} \left( \| u^n \|_{H^1(F^n(t))}^2 \: + \:  \| u^n \|^2_{L^2(S^n(t))} \right) \, dt \: \le \:  C  
\end{array}\label{unn2}  
\end{multline}
yields that 
\begin{equation}
 \int_0^T \| (u^n - P^n_S u^n)(t,\cdot) \cdot \nu \|_{L^p(\pa S^n(t))}^2 \, dt \:  \rightarrow  \: 0, \quad \forall \, p < 4. 
 \end{equation}
Eventually, we get that 
\begin{equation} \label{boundvn}
  \| (1- \chi^n_S) (v^n -  u^n) \|_{L^2(0,T; L^p(\Omega))} \rightarrow 0, \quad \forall \, p < 4. 
\end{equation}
This will be much important in the treatment of the nonlinear terms. 

\medskip
We conclude this paragraph with the construction of the operator $V^{\delta,n}$ satisfying \eqref{extra}. We take $U,U_S$ in $H^1_\sigma(\R^3\setminus S_0) \times H^1_\sigma(\overline{S_0})$. Up to  an  extension of $U_S$, there is no loss of generality assuming that $U_S \in H^1_\sigma(\R^3)$.

\medskip
{\em Step 1.} We shall construct a field $V$ such that $\div V$ = 0, 
\begin{equation} \label{BCW}
 V\vert_{\pa S_0}  = U_S  + (U - U_S) \cdot \nu \, \nu, \quad \mbox{ and } \: V\vert_{\pa (S_0)_\delta} =  U.  
\end{equation}
Therefore, we introduce a system of orthogonal curvilinear coordinates $(s_1,s_2,z)$ in a tubular neighborhood of $\pa S_0$: $s_1, s_2$ are coordinates along the surface $\pa S_0$, whereas $z$ denotes a transverse coordinate. In particular, $\pa S_0 = \{ z = 0 \}$.  We set 
$$ e_1 \: := \: \frac{1}{h_1} \frac{\pa}{\pa s_1}, \quad    e_2 \: := \:  \frac{1}{h_2} \frac{\pa}{\pa s_2}, \quad e_z  :=  \nu \: = \:  \frac{1}{h_z} \frac{\pa}{\pa z} $$
the associated orthonormal vectors, with scale factors $h_1, h_2, h_z \ge 0$. We remind that 
\begin{equation} \label{gradient}
 \na f = \frac{1}{h_1} \pa_{s_1} f e_1 \: + \: \frac{1}{h_2} \pa_{s_2} f e_2  \: + \: \frac{1}{h_z} \pa_{z} f e_z  
 \end{equation}
for a scalar function $f$, whereas 
\begin{equation} \label{divergence}
 \div f = \frac{1}{h_1 h_2 h_z} \left( \pa_{s_1} (h_2 h_z f_1) \: + \:  \pa_{s_2} (h_1 h_z f_2) \: + \:  \pa_z (h_1 h_2 f_z) \right)
 \end{equation}
for any field $f = f_1 e_1 + f_2 e_2 + f_z e_z$. 
We then set 
$$ V_1 \: := \: \left( 1 - \chi(n z) \right) U + \chi(n z) (U_S + [(U- U_S) \cdot e_z] \, e_z ) $$
for  a smooth  truncation function $\chi : \mathbb R \to [0,1]$ equal to $1$ in a neighborhood of $0$. Clearly, for all  $p \le  6$, and $\frac{1}{q} + \frac{1}{6} = \frac{1}{p}$, 
\begin{equation} \label{estimW1}
\begin{aligned}
 \| V_1 - U \|_{L^p((S_0)_\delta \setminus S_0)} \:  & \le \: C_{p, \delta} \,  n^{-1/q} \,  \| (U,U_S) \|_{L^{6}((S_0)_\delta \setminus S_0)} \\ 
 & \le   \, C'_{p,\delta} \,  n^{-1/q} \,  \| (U,U_S) \|_{H^1((S_0)_\delta \setminus S_0)}. 
\end{aligned}
\end{equation}
Also, $V_1$ satisfies \eqref{BCW}. But it is not divergence-free: formula \eqref{divergence} yields   
$$  \div V_1 \: = \: \chi(n z) \, \div([(U- U_S) \cdot e_z] e_z ) $$
so that for all $p \le 2$, 
$$  \| \div V_1 \|_{L^p((S_0)_\delta \setminus S_0)} \: \le \: C \ n^{1/2-1/p} \| U - U_S \|_{H^1((S_0)_\delta \setminus S_0)}. $$
To obtain a divergence-free field,  we note that both $U$ and $U_S$ have zero flux through $\partial S_0$ and  \cite[Theorem 3.1]{Galdibooknew}:  there exists a field $V_2$ such that 
$$ \div V_2 = - \div V_1 \quad  \mbox{ in } \quad (S_0)_\delta \setminus S_0, \quad V_2 \vert_{\pa S_0} = V_2\vert_{\pa (S_0)_\delta} = 0, $$
and for all $p \in ]1,2]$,
$$ \| V_2 \|_{W^{1,p}((S_0)_\delta \setminus S_0)}\: \le \: C_\delta \ n^{1/2-1/p} \, \| U - U_S  \|_{H^1((S_0)_\delta \setminus S_0)}. $$ 
In particular, by Sobolev imbedding, one has for all $p_* \le  6$ 
\begin{equation} \label{estimW2}
 \| V_2  \|_{L^{p_*}((S_0)_\delta \setminus S_0)}\: \le \: C_\delta \, n^{1/6-1/p_*}  \, \| (U,U_S)  \|_{H^1((S_0)_\delta \setminus S_0)}. 
\end{equation}
Finally, the field $\: V \: := \:  V_1 + V_2$ fulfills our requirements.  

\medskip
{\em Step 2.} We construct a field $W$ such that $\div W$ = 0, 
\begin{equation} \label{BCtildeW}
W\vert_{\pa S_0}  = [(U - U_S) \cdot \nu] \, \nu, \quad \mbox{ and } \: W\vert_{\pa (S_0)_\delta} =  0.  
\end{equation}
In the same spirit as in the first step, we take 
$$ W_1 \: := \: \chi(\frac{2 z}{\delta})   [\left( U - U_S \right) \cdot \nu\vert_{z=0}] \, e_z $$
where $\chi$ is again a truncation function: $\chi = 1$ near $0$, and $\chi = 0$ outside $[-1,1]$. 
A rapid computation shows that 
\begin{equation} \label{estimW3}
 \| W_1 \|_{L^p((S_0)_\delta \setminus S_0)} \: \le \: C_\delta      \| (U- U_S) \cdot \nu  \|_{L^p(\pa S_0)}, \quad \forall \, p 
 \end{equation}
By Proposition \ref{Bogovski}, there exists a field $W_2$ such that 
$$ \div W_2 = - \div W_1 \quad  \mbox{ in } \quad (S_0)_\delta \setminus S_0, \quad  W_2 \vert_{\pa S_0} = W_2\vert_{\pa (S_0)_\delta} = 0, $$
and
$$ \| W_2  \|_{H^1((S_0)_\delta \setminus S_0)}\: \le \: C_\delta  \, \| (U - U_S) \cdot \nu   \|_{L^2(\pa S_0)}. $$ 
In particular, by Sobolev imbedding, one has for all $p \le  6$ 
\begin{equation} \label{estimW4}
 \| W_2  \|_{L^p((S_0)_\delta \setminus S_0)} \: \le \: C_\delta  \, \| (U - U_S) \cdot \nu   \|_{L^2(\pa S_0)}. 
\end{equation}
Finally, the field $W \: := \:  W_1 + W_2$ fulfills our requirements.  

\bigskip
Eventually, we set 
\begin{equation*}
\left\{
\begin{aligned}
V^\delta[U,U_S] \: & = \: U \:  \mbox{ outside }  \:  (S_0)_\delta, \\
V^\delta[U,U_S] \: & = \: V - W \:  \mbox{ in }  \: (S_0)_\delta \setminus S_0, \\
V^\delta[U,U_S] \: & = U_S  \:  \mbox{ in }  \: S_0.
\end{aligned}
\right.
\end{equation*}
Combining \eqref{estimW1}, \eqref{estimW2}, \eqref{estimW3} and \eqref{estimW4} leads to \eqref{extra}. 

\subsection{Approximation of the test functions}
The weak formulation of the momentum equation involves discontinuous test functions $\varphi \in {\cal T}_T$:
$$ \varphi = (1- \chi_S) \varphi_F + \chi_S \varphi_S, \quad \varphi_F \in {\cal D}([0,T);  {\cal D}_\sigma(\overline{\Omega})), \quad  \varphi_S \in {\cal D}([0,T);  {\cal R}), $$
with 
$$ \varphi_F \cdot \nu\vert_{\pa \Omega} = 0, \quad  \varphi_F \cdot \nu\vert_{\pa S(t)} = \varphi_S \cdot \nu\vert_{\pa S(t)} \quad \forall \, t.$$ 
On the contrary, the approximate momentum equation {\bf c)} involves continuous (or at least $H^1$) test functions. Hence, we will have to approach $\varphi $ by a sequence $\displaystyle (\varphi^n)$ in   ${\alert L^2}(0,T; H^1_\sigma(\Omega))$. Due to the discontinuity of the limit, the $\displaystyle \varphi^n(t,\cdot)$'s will exhibit strong gradients near $\displaystyle \pa S^n(t)$. Precise estimates are needed, that  are the purpose of 
\begin{proposition} \label{approximatetest}
Let $\alpha > 0$. There exists a sequence $(\varphi^n)$  in  
$W^{1,\infty}(0,T ; L^2_{\sigma}(\Omega)) \cap L^{\infty}(0,T; H^1_\sigma(\Omega))$, of the form
$$ \varphi^n \:  = \:  (1 - \chi^n_S) \varphi_F  \: + \:  \chi^n_S \varphi^n_S, $$  
that satisfies 
\begin{itemize}
\item $\| \sqrt{\chi^n_S} (\varphi^n_S - \varphi_S) \|_{C([0,T]; L^p(\Omega))} = O(n^{-\alpha/p})$ for all $p \in [2,6]$.  
\item $\varphi^n \rightarrow \varphi \: \mbox{ strongly in } \: C([0,T];  L^6(\Omega))$. 
\item $\| \varphi^n \|_{C([0,T]; H^1(\Omega))} = O(n^{\alpha/2})$. 
\item $\|  \chi^n_S (\pa_t + P^n_S u^n \cdot \na) \left(\varphi^n -   \varphi_S \right) \|_{ L^{\infty}(0,T; L^6(\Omega))} = O(n^{-\alpha/6})$. 
\item $ (\pa_t + P^n_S u^n \cdot \na) \varphi^n \: \rightarrow \: (\pa_t + P_S u \cdot \na) \varphi \: \mbox{ weakly * in } \quad L^\infty(0,T; L^6(\Omega))$. 
\end{itemize}
\end{proposition}

\medskip
{\em Proof of the proposition}. 

\medskip
The point is to build a good approximation $\varphi^n_S$ of $\varphi_S$ over the solid domain. Broadly, we  want
$$ \varphi^n_S(t, \cdot)\vert_{\pa S^n(t)} = \varphi_F(t,\cdot)\vert_{\pa S^n(t)}   \quad \forall \, t, $$
and 
$$ \varphi^n_S(t, \cdot) \approx \varphi_S(t,\cdot) \quad  \mbox{ in $S^n(t)$  away from  a $n^{-\alpha}$ neighborhood of } \pa S^n(t) \quad \forall \, t.  $$
Therefore, we proceed as for $v^n$, by using lagrangian coordinates: we define $\Phi_S$ and $\Phi_F$ through the formulas 
$$ \varphi_S\left(t,\phi^n_{t,0}(y)\right) \: = \: d\phi^n_{t,0}\vert_{y}\left(\Phi_S(t,y) \right) , \quad \varphi_F\left(t,\phi^n_{t,0}(y)\right) \: = \: d\phi^n_{t,0}\vert_{y}\left(\Phi_F(t,y) \right),  $$
and the goal is to define properly some $\Phi^n_S$, related to $\varphi^n_S$ by the formula 
$$  \varphi^n_S\left(t,\phi^n_{t,0}(y)\right) \: = \: d\phi^n_{t,0}\vert_{y}\left(\Phi^n_S(t,y) \right). $$
{\em Note that $\Phi_S$ and $\Phi_F$ depend on $n$ through the propagator $\phi^n$, but we omit it from our notations.} The only thing we have to keep in mind is that the bounds on $\phi^n$ guarantee that 
$\Phi_S$ and $\Phi_F$ are uniformly  bounded in $\displaystyle W^{1,\infty}(0,T; H^k_{loc}(\R^3))$ for all $k.$ 

\medskip
Thanks to the change of coordinates, the problem  is now in the fixed domain $S_0$.  Roughly, we  want to build $\Phi^n_S$ in such a way that  
$$ \Phi^n_S(t, \cdot)\vert_{\pa S_0} = \Phi_F(t,\cdot)\vert_{\pa S_0}   \quad \forall \, t, $$
and 
$$ \Phi^n_S(t, \cdot) \approx \Phi_S(t,\cdot) \quad  \mbox{ in $S_0$  away from  a $n^{-\alpha}$ neighborhood of } \pa S_0 \quad \forall \, t.  $$
Note that time is only a parameter in the system. The construction of $\Phi^n_S$  follows the one of $V$, performed in the previous paragraph, Step 1. We take $\Phi^n_S$ under the form 
$$ \Phi^n_S \: = \:  \Phi^n_{S,1} \: + \:   \Phi^n_{S,2}. $$
The first term has the explicit form 
$$ \Phi^n_{S,1} \: = \:  \Phi_S \: + \: \chi(n^\alpha \, z) \left( (\Phi_F - \Phi_S) - [(\Phi_S - \Phi_F) \cdot e_z] \, e_z  \right).  $$  
Again, $\chi$ is a smooth truncation function near $0$, and $z$ is a coordinate transverse to the boundary: $\pa S_0 = \{ z = 0 \}$.  
It is easily seen that $\Phi^n_{S,1}$ satisfies the right boundary condition at $\pa S_0$. Moreover, 
\begin{equation} \label{estimpsi1}
 \| \Phi^n_{S,1} - \Phi_S \|_{W^{1,\infty}(0,T; L^p(S_0))} \: \le \: C \, n^{-\alpha/p} \quad \forall \, p < \infty,  \quad  \| \Phi^n_{S,1} - \Phi_S \|_{W^{1,\infty}(0,T; H^1(S_0))} \: \le \: C \,  n^{\alpha/2}.
\end{equation}
But it is not divergence-free. By applying formula \eqref{divergence}, we get 
$$ \div  \Phi^n_{S,1} \: = \:  \chi(n^\alpha \, z) \j^n, \quad \j^n \: := \:   \div \left( x \mapsto  \left( (\Phi_F - \Phi_S) - [(\Phi_S - \Phi_F) \cdot e_z] \, e_z  \right)\right). $$
In particular,  $\j^n$  is uniformly bounded in $W^{1,\infty}(0,T; L^2(S_0))$  

\medskip
By Proposition \ref{Bogovski}, there exists some field $ \Phi^n_{S,2}$ satisfying 
$$ \div \Phi^n_{S,2} = - \div \Phi^n_{S,1} \quad \mbox{ in } \: S_0, \quad  \Phi^n_{S,2}\vert_{\pa S_0} = 0, $$
and 
\begin{equation} \label{estimpsi2}
\| \Phi^n_{S,2} \|_{W^{1,\infty}(0,T; H^1(S_0))} \: \le  \: C \, \| \chi(n^\alpha z) \, \j^n \|_{W^{1,\infty}(0,T; L^2(S_0))} \: \le \: C \, n^{-\alpha/2}. 
\end{equation} 
In particular,
\begin{equation} \label{estimpsi3}
\| \Phi^n_{S,2} \|_{W^{1,\infty}(0,T; L^6(S_0))} \: \le \: C \, n^{-\alpha/2}.
\end{equation}  

\medskip
Back to the moving domain (in variable $x$), one can combine the estimates \eqref{estimpsi1}-\eqref{estimpsi2}-\eqref{estimpsi3} with the uniform bound on $\phi^n$ in $W^{1,\infty}(0,T; \, C^{\infty}(\Omega))$. From there, one can deduce the estimates of the proposition. For the sake of brevity, we only treat the two last items.  Namely, we write 
\begin{align*}
 \| \chi^n_S (\pa_t + P^n_S u^n \cdot \na) \left(\varphi^n -   \varphi_S \right) \|_{L^{\infty}((0,T); L^6(\Omega))}  &\: \le \:  C \, \| \frac{\pa}{\pa t}  \, d\phi^n_{t,0}\vert_{y}\left( \Phi^n_S - \Phi_S \right) \|_{L^{\infty} ((0,T); L^6(S_0))} \\
& \: \le \: C \, n^{-\alpha/6}\,,
\end{align*}
where the last inequality involves \eqref{estimpsi1} and \eqref{estimpsi3}. This bound implies in turn that 
\begin{equation*}
 ( \pa_t + P^n_S u^n  \cdot \na) \varphi^n \: 
   = \:  (1 - \chi^n_S) ( \pa_t + P^n_S u^n  \cdot \na) \varphi_F \: + \:  \chi^n_S  ( \pa_t + P^n_S u^n \cdot \na) \varphi_S \: + \: O(n^{-\alpha/6}) \quad \mbox{ in } \:  L^6(\Omega) 
\end{equation*}
The products at the r.h.s. are then easily handled using the strong convergence of $\chi^n_S$ (and the weak convergence of $P^n_S u^n$). 
We obtain 
$$ ( \pa_t + P^n_S u^n \cdot \na)  \varphi^n  \: \rightarrow (\pa_t + P_S u \cdot \na ) \varphi  \quad \mbox{weakly * in } \:  L^\infty(0,T; L^6(\Omega))$$
as expected. This concludes the proof of the proposition.

\subsection{Convergence in the momentum equation: linear terms}
We now have all the elements to study the asymptotics of the approximate momentum equation {\bf c)}.  Given an arbitrary $\varphi \in {\cal T}_T$, we consider an approximate sequence $(\varphi^n)$ as in Proposition  \ref{approximatetest}.  We shall  take $\varphi^n$ as a test function in {\bf c)}, and let $n$ tend to infinity, so as  to recover \eqref{momentum}. We shall rely on the fields $u^n_F$ and  $u_F$ introduced in paragraph \ref{aprioribounds}. We remind that 
\begin{equation} \label{propertiesuf}
(1 - \chi^n_S ) \, u^n_F \: = \:  (1 - \chi^n_S ) \, u^n, \quad u^n_F \: \rightarrow \: u_F \quad \mbox{ weakly in } \:  L^2(0,T; H^1_\sigma(\Omega)). 
\end{equation}
To lighten notations, we shall write  $u^n_S := P^n_S u^n$, $ \: u_S := P_S u$.  We remind that these rigid fields satisfy 
\begin{equation} \label{propertiesus}
 u^n_S  \: \rightarrow \: u_S \quad \mbox{ weakly * in } \:  L^\infty(0,T; W^k_{loc}(\R^3)) \quad \forall k.
 \end{equation}

\medskip
 In this paragraph, we consider the asymptotics of all terms but the convection one.
\begin{itemize}
\item We write the diffusion term as: 
\begin{align*}
& \int_0^T \int_\Omega 2 \mu^n D(u^n) : D(\varphi^n)   = \int_0^T \int_{\Omega} \left( 2 \mu_F (1 - \chi^n_S)  D(u^n_F) + \frac{1}{n^2}\chi^n_S D(u^n) \right) : D(\varphi^n) \\
&  =  \int_0^T  \int_{\Omega} 2 \mu_F  (1 - \chi^n_S)  D(u^n_F)  : D(\varphi_F) \: + \: \frac{1}{n^2} \int_0^T \int_\Omega \chi^n_S \, D(u^n_S) : D(\varphi^n) \: := \: I_1^n \: + \: I_2^n. 
 \end{align*}
 From the strong convergence of $\chi^n_S$ to $\chi_S$ in $C([0,T]; L^p(\Omega))$, and the weak convergence of $u^n_F$ to $u_F$ in $L^2(0,T; H^1(\Omega))$, we deduce 
 $$ I_1^n \: \rightarrow \:  \int_0^T  \int_{\Omega} 2 \mu_F  (1 - \chi_S)  D(u_F)  : D(\varphi_F). $$
As regards $I_2^n$, we  use the bounds 
$$  \| \sqrt{\mu^n} D(u^n) \|^2_{L^2((0,T) \times \Omega)} = O(1), \quad     \| \varphi^n \|_{L^{\infty}(0,T ; H^1(\Omega))} = O(n^{\alpha/2})  $$
established in the previous paragraphs. They imply
$$ | I_2^n | \: \le \: \frac{C}{n^2} \,  \| \chi_S D(u^n) \|_{L^2((0,T) \times\Omega )}  \, \| D(\varphi^n) \|_{L^2((0,T) \times \Omega)} \: \le \: \frac{C}{n^{1-\alpha/2}}  $$
If we choose $\alpha  < 2$, then $I_2^n$ goes to zero as $n$ goes to infinity, and finally
$$ \int_0^T \int_\Omega 2 \mu^n D(u^n) : D(\varphi^n)  \: \rightarrow \:\, \int_0^T  \int_{F(t)}    2 \mu_F  \, D(u_F)  : D(\varphi_F). $$
\item  The boundary term at $\pa \Omega$ reads 
\begin{align*}
  \frac{1}{2\beta_\Omega}\int_0^T \int_{\pa \Omega}  (u^n \times \nu) \cdot (\varphi^n \times \nu) & =  \frac{1}{2\beta_\Omega}\int_0^T \int_{\pa \Omega}  (u^n_F \times \nu) \cdot (\varphi_F \times \nu) \\
& \rightarrow \frac{1}{2\beta_\Omega}\int_0^T \int_{\pa \Omega}  (u_F \times \nu) \cdot (\varphi_F \times \nu)
\end{align*}
by the weak convergence of $u^n_F$ in $L^2(0,T ; H^1(\Omega))$. 
\item We deal with the boundary term at $\pa S^n$ as in the Galerkin approximation. 
We introduce $u^n_S := P^n_S u^n,$ $r^n_S := P^n_S \varphi^n = P^n_S  \varphi^n_S $ and capital letters to denote velocity fields when seen through the change of variable. 
We then have, as in the computation for the Galerkin method : 
\begin{align*}
& \frac{1}{2\beta_S} \int_0^T  \int_{\pa S^n(t)}  ((u_F^n -  u^n_S) \times \nu) \cdot ((\varphi_F^n -  r^n_S) \times \nu)  \\
=  & \frac{1}{2\beta_S} \int_0^T  \int_{\pa S_0}  ((U^n_F - U^n_S) \times \nu) \cdot ((\Phi_F -  R^n_S) \times \nu   )\,,
\end{align*}
where we used that $\varphi_F^n = \varphi_F.$
We note here that $\varphi^n$ converges to $\varphi$ in $C([0,T];L^6(\Omega))$ so that combining with the strong convergence
of $\chi^n_S$ it yields that $r^n_S$ converges to $r_S := P_S\varphi$ in $L^2(0,T; H^1_{loc}(\mathbb R^3)).$
Through the change of variable, Lemma \ref{lemmawn} yields that :
$$
R^n_S \to R_S \text{ strongly in } L^2(0,T;H^{1/2}(\partial S_0)).
$$
Then, we combine the respective convergences of $u^n_F,u^n_S$  with Lemma \ref{lemmawn}
yielding,  with obvious notations : 
$$
U^n_F \to  U_F\text{ weakly in } L^2(0,T;H^1(\Omega)),\quad
U^n_S \to  U_S \text{ weakly in } L^2(0,T;H^1(\Omega))\,.
$$
We apply these convergences together with the continuity of traces on $\partial S_0 \subset \subset \Omega,$ and go back to the moving geometry, to obtain finally :
\begin{multline*}
\frac{1}{2\beta_S} \int_0^T  \int_{\pa S^n(t)}  ((u^n -  u^n_S) \times \nu) \cdot ((\varphi^n -  \varphi^n_S) \times \nu) \\
\to 
 \frac{1}{2\beta_S} \int_0^T  \int_{\pa S(t)}  ((u_F -  u_S) \times \nu) \cdot ((\varphi_F -  \varphi_S) \times \nu) \,.
\end{multline*}
\item to treat the penalization term we use the bounds 
$$   n  \|  \sqrt{\chi^n_S} (u^n - P_S^n u^n) \|^2_{L^2((0,T) \times \Omega)} = O(1), \quad  \| \sqrt{\chi^n_S} (\varphi^n_S - \varphi_S) \|_{C([0,T]; L^6(\Omega))} = 
O(n^{-\alpha/2}) $$
established in the previous paragraph.  We also remind that  $\varphi_S$, as a rigid vector field, satisfies 
 $\varphi_S = P^n_S \varphi_S$.  From there, 
\begin{align*}
& \left| n \int_0^T \int_\Omega  \chi^n_S  (u^n - P^n_S u^n)  \cdot (\varphi^n - P^n_S \varphi^n) \right| \: \\
=  \: & 
 \left| n \int_0^T \int_\Omega  \chi^n_S  (u^n - P^n_S u^n)  \cdot \left((\varphi^n_S - \varphi_S) - P^n_S (\varphi^n - \varphi_S)\right) \right| \\
 = \: & \left| n \int_0^T \int_\Omega  \chi^n_S  (u^n - P^n_S u^n)  \cdot (\varphi^n_S - \varphi_S) \right| \: \le \: C n^{1/2-\alpha/2} 
 \end{align*}
 If we choose $\alpha > 1$ (which is compatible with the former constraint $\alpha < 2$), the penalization term vanishes as $n \rightarrow +\infty$. 
\end{itemize}

\subsection{Strong convergence of $(u^n)$}

To show that $(S,u)$ is a weak solution over $(0,T)$, we still have to pass to the limit in the convection term
$$ {\rm conv}^n \: := \:  \int_0^T \int_\Omega \rho^n \, \left(u^n \cdot \pa_t \varphi^n + v^n \otimes u^n : \na \varphi^n \right). $$
To compute this limit, we first prove 
\begin{proposition} \label{prop_strongcvgce}
Up to the extraction of a subsequence,  $(u^n)$ converges to $u$ in $L^2((0,T) \times \Omega).$
\end{proposition}
 This result is obtained applying the method introduced in the reference \cite{SanMartin&al02} (see also \cite{Feireisl03} for the 3D case). 
We first introduce  some notations. Given $0 \le s \leq 1$ and $S$ a bounded connected subset
 $\Subset \Omega$, we denote 
 $$ \mathcal{R}^s[S] \: = \: \mbox{ 
the closure of } \:  \{ v \in H^1_{\sigma}(\Omega) \: \text{ such that } v_{|_{S}}  \in \mathcal{R}\} \,  
\: \text{in}  \: H^s(\Omega).$$ 
As $\mathcal{R}^s[S]$  is a closed subspace of $H^s(\Omega)$ we denote 
$P^{s}[S]$  the orthogonal projector from $H^s(\Omega)$ onto this subspace. Given $s'>s,$ we recall that   
$\mathcal{R}^{s'}[S]$ is a dense subspace of  $\mathcal{R}^{s}[S]\,,$ and that the imbbeding 
$\mathcal{R}^{s'}[S] \subset \mathcal{R}^{s}[S]$ is  compact. If $s =0,$ we shall drop exponent $s$.
We emphasize that in the case $s=0$ the projector $P[S]$ does not coincide with the $P_S$ introduced in \eqref{defPs}.  

\medskip

Our first step is the following approximation lemma :
\begin{lemma}  \label{lemvhbar} Let $s < \frac{1}{3}$.
\begin{itemize} 
\item[\bf i)]  The sequence  $(u_n)$ is uniformly bounded in $L^2(0,T;H^s(\Omega))$. Moreover,  there is $\eps = \eps(s) > 0$ such that  for all $h < \delta/2$,
\begin{equation} \label{eq_hbarun}
\int_0^T \|u^n(t,\cdot) - P^s[(S^n(t))_h] \,  u^n(t,\cdot)\|^2_{H^s(\Omega)} \: \leq  \: C \left(h^\eps \, + \, n^{-\eps} \right). 
\end{equation}
\item[\bf ii)] One has $u \in L^2(0,T;H^s(\Omega))$. Moreover,   there exists $\eps = \eps(s)$ such that  for all $h < \delta/2$,
\begin{equation} \label{eq_hbaru}
\int_0^T \|u(t,\cdot) - P^s[(S(t))_h] \,  u(t,\cdot)\|^2_{H^s(\Omega)}\text{d$t$} \: \leq \:  C h^\eps, 
\end{equation}
\end{itemize}
where, in both cases, the constant $C$ depends only on initial data.
\end{lemma}
{\em Proof of the lemma.}
We only prove the first item of the lemma, the second one being simpler. It relies on  the construction of a suitable approximation $v^n_h$ of $u^n$, rigid in a $h$-neighborhood of $S^n$. This approximation will satisfy the following properties :
\begin{itemize}
\item $(v^n_{h})$ is bounded in  $L^2(0,T;H^s(\Omega))$ for $s$ small  enough  
\item  $v^n_{h}(t,\cdot) = P_{S^n} u^n(t,\cdot)$ in $(S^n(t))_{h}$ and $v^n_{h}(t,\cdot) = u^n(t,\cdot)$ outside $(S(t))_{\delta}$ for a.a. $t \in (0,T)$.

\smallskip
 Note  that it implies  $v^n_h(t,\cdot) \in \mathcal R^s[(S^n(t))_h]$ for a.a. $t \in (0,T)$. 
\item for $h$ sufficiently small and for a.a. $t\in (0,T)$ there holds
\begin{equation} \label{u-vh}
\begin{aligned}
 \|u^n(t,\cdot) - v^n_{h}(t,\cdot)\|_{L^2(\Omega \setminus (S^n(t))_h)}& \: \le \:  C \, h^{\frac{1}{3}} \left( \|P^n_S \, u^n(t,\cdot)\|_{L^2(\Omega)} + \|u^n(t,\cdot)\|_{H^1(F^n(t))}\right)   \\
 & \quad + C \, \| (u^n - P^n_S u^n) \cdot \nu \|_{L^2(\pa S^n(t))}\,, \\
\vspace{1cm} \|\ v^n_{h}(t,\cdot)\|_{H^s(\Omega)}  & \: \le \:    C  (1 +  h^{\frac{1}{3}-s}) \left( \|P^n_S \, u^n(t,\cdot)\|_{L^2(\Omega)} + \|u^n(t,\cdot)\|_{H^1(F^n(t))}\right)   \\
 & \quad +  C \, \| (u^n - P^n_S u^n) \cdot \nu \|_{L^2(\pa S^n(t))}.
\end{aligned} 
\end{equation}
\end{itemize}
Before giving  further details  on the construction of $v^n_h$ we explain how the previous properties
imply Lemma \ref{eq_hbarun}. By interpolation of \eqref{unn1} and \eqref{unn2}, we obtain 
\begin{equation} \label{unn3}
 \int_0^T \| (u^n - P^n_S u^n)(t,\cdot) \cdot \nu \|_{L^2(\pa S^n(t))}^2 \, dt \: \le \: \frac{C}{\sqrt{n}}. 
 \end{equation}
We square the inequalities  in \eqref{u-vh} and integrate from $0$ to $T$. Using \eqref{unn3} with the uniform bounds  
 \eqref{boundun}, we end up with 
\begin{equation} \label{u-vh2}
 \left(  \int_0^T \| u^n - v^n_h \|^2_{L^2(\Omega \setminus (S^n(t))_h)}\right)^{1/2} \: \le \:  C \left( h^{1/3} + \frac{1}{\sqrt{n}} \right), \quad   \| v^n_h \|_{L^2(0,T;H^s(\Omega))} \: \le \: C, \quad \forall s \le 1/3.  
 \end{equation}
 Moreover, 
 \begin{align*}
&  \left(  \int_0^T  \| u^n - v^n_h \|^2_{L^2((S^n(t))_h)}\right)^{1/2}    = \:  \left(  \int_0^T \| u^n - P^n_S u^n \|^2_{L^2((S^n(t))_h)}\right)^{1/2} \\
& \le \: \left(  \int_0^T \| u^n - P^n_S u^n \|^2_{L^2(S^n(t))}\right)^{1/2} \: + \:  \left(  \int_0^T \| u^n - P^n_S u^n \|^2_{L^2((S^n(t))_h \setminus S^n(t))}\right)^{1/2}
  \end{align*}
Using \eqref{boundun}, we get 
 \begin{align*}
  \left(  \int_0^T \| u^n - v^n_h \|^2_{L^2((S^n(t))_h)}\right)^{1/2} 
& \le \:  \frac{C}{\sqrt{n}} \: + \: \left(  \int_0^T \| u^n - P^n_S u^n \|^2_{L^2((S^n(t))_h \setminus S^n(t))}\right)^{1/2} \\
& \le \:  \frac{C}{\sqrt{n}} \: + \:  C \sqrt{h} \, \left( \int_0^T  \| u^n - P^n_S u^n \|^2_{L^4((S^n(t))_h \setminus S^n(t))} \right)^{1/2}  \\
& \le \: \frac{C}{\sqrt{n}} \: + \:  C \sqrt{h} \, \left( \int_0^T  \| u^n - P^n_S u^n \|^2_{H^1(F^n(t))} \right)^{1/2} \: \le \: 
 \frac{C}{\sqrt{n}} \: + \:  C \sqrt{h} \, 
 \end{align*}
 Combining this last inequality with the first inequality in \eqref{u-vh2} yields 
 \begin{equation} \label{u-vh3}
  \| u^n - v^n_h \|_{L^2((0,T) \times \Omega)} \: \le \: C \, \left( h^{1/3} + \frac{1}{\sqrt{n}} \right)
  \end{equation} 
 As regards the $H^s$ norm, $s \le  1/3$, we use the second inequality in \eqref{u-vh2} to write 
 \begin{align*}
  \| u^n - v^n_h \|_{L^2(0,T; H^s(\Omega))} \: & \le \: \left( \int_0^T \| u^n - P^n_S u^n \|_{H^s(S^n(t))}^2 dt \right)^{1/2} 
  +   \left( \int_0^T \| u^n \|_{H^s(F^n(t))}^2 dt \right)^{1/2} \\
  & \quad + \:  \left( \int_0^T \| v^n_h \|_{H^s(F^n(t))}^2 dt \right)^{1/2} \: \le \: \| u^n - P_S^n u^n \|_{L^2(0,T; H^s(S^n))} \: + \: O(1) 
 \end{align*}
Finally, we have 
\begin{align*} 
\| u^n - P_S^n u^n \|_{L^2(0,T; H^s(S^n))}  & \: \le \: C \,
\| u^n - P^n_S u^n  \|_{L^2((0,T) \times S^n)}^{1 - s} \,   \| u^n - P^n_S u^n  \|_{L^2(0,T; H^1(S^n))}^s  \\
& \: \le \: C \, \left( \frac{1}{\sqrt{n}} \right)^{1-s} \, n^s \: \le \:  C  \quad  \mbox{as soon as } \:  s \le  \frac{1}{3}.  
&
\phantom{toto}
\end{align*}
We end up with 
\begin{equation} \label{u-vh4}
\| u^n - v^n_h \|_{L^2(0,T;H^s(\Omega))} \: \le \: C  \quad  \mbox{ as soon as } s \le  \frac{1}{3}.  
\end{equation}
One last interpolation between \eqref{u-vh3}  and \eqref{u-vh4} shows that for all $s < 1/3$ and $\eps = \eps(s) > 0$, 
$$ \| u^n - v^n_h \|_{L^2(0,T; H^s(\Omega))} \: \le \: C \left( h^{\eps} + n^{-\eps} \right). $$
 As $v^n_h(t, \cdot)$ belongs to ${\cal R}^s[(S^n(t))_h]$ for all $t$, by definition of the projection, the same inequality holds replacing $v^n_h$ by $P[(S^n(t))_h]$, as expected.

\bigskip
\bigskip
We still have to achieve the  construction of $v^n_{h}$. It follows the  construction of $v^n$, {\it cf } paragraph \ref{aprioriboundsv}.  It is actually simpler, because we only look for a $v^n_h$ with $H^s$ regularity for small $s$. In particular, jump on the tangential part at $\pa (S^n)_h$ and $\pa (S^n)_\delta$ will be allowed. 

\medskip
As before,  we  go back to Lagrangian coordinates : we look for a $v^n_h$ under the form 
$$
 v^n_{h}(t,\phi_{t,0}(y)) = d\phi_{t,0}|_{y} V^n_{h}(t,y).
$$
Also, we define $U^n$ and $U^n_S$ through 
 $$ E_\Omega u^n\left(t,\phi^n_{t,0}(y)\right) \: = \: d\phi^n_{t,0}\vert_{y} \, U^n(t,y), \quad P^n_S u^n\left(t,\phi^n_{t,0}(y)\right) \: = \: d\phi^n_{t,0}\vert_{y} \, U^n_S(t,y).  $$
In this way, we are back to a static problem. {\em For brevity, we shall omit temporarily the time dependence in our  notations}. The point is to  build a field $V^n_h$ satisfying 
$$ V^n_h = U^n_S \quad  \mbox{ in } \:  (S_0)_h, \quad V^n_h = U^n \quad  \mbox{ outside } (S_0)_\delta,  $$
and suitable estimates. 

\medskip
Therefore, we follow paragraph \ref{aprioriboundsv}. 
We parametrize $\displaystyle (S_0)_{\delta} \setminus S_0$  by curvilinear coordinates $(s_1,s_2,z)$,  $z$ being the distance at $\partial S_0$.  Hence, 
$\displaystyle \partial (S_0)_{h} = \{  z= h \}.$  Then, we introduce 
$$ V^n_{h,1}  \: := \: \left( 1 - \chi\left(\frac{z-h}{h}\right) \right) U^n + \chi\left(\frac{z-h}{h}\right) (U^n_S + [(U^n- U^n_S) \cdot e_z] \, e_z ) $$
and the solution  $V^n_{h,2}$ of 
$$
\left\{
\begin{array}{rcll}
{\rm div} V^n_{h,2} &=& - {\rm div} V^n_{h,1}\,, & \text{ in $(S_0)_{\delta} \setminus (S_0)_{h}$}\,, \\[4pt]
V^n_{h,2}& =& 0 & \text{ on $\partial (S_0)_{\delta}$ and $\partial (S_0)_{h}$}  
\end{array}
\right.
$$
Computations similar to those of paragraph \ref{aprioriboundsv}  yield :
\begin{eqnarray} \label{eq_boundV10}
\|V^n_{h,1} - U^n\|_{L^2((S_0)_{\delta} \setminus (S_0)_{h})} &\leq&  h^{\frac{1}{3}} \|(U^n,U^n_S)\|_{H^1((S_0)_{\delta} \setminus (S_0)_{h}) \times {\cal R}},  \\
\|V^n_{h,1} \|_{H^1((S_0)_{\delta} \setminus (S_0)_{h})} & \leq &  h^{-\frac{2}{3}}  \|(U^n,U^n_S)\|_{H^1((S_0)_{\delta} \setminus (S_0)_{h})\times {\cal R}}, 
\label{eq_boundV11} \\
\|V^n_{h,2}\|_{H^1((S_0)_{\delta} \setminus (S_0)_{h})} &\leq&  C \ \|(U^n,U^n_S)\|_{H^1((S_0)_{\delta} \setminus (S_0)_{h})\times {\cal R}}.
\label{eq_boundV20}
\end{eqnarray}
Let us emphasize that the constant $C$ in the last inequality can be chosen uniformly in $h$, see \cite[Theorem III.3.1]{Galdibooknew}. 
It follows by interpolation that 
\begin{equation} \label{interpol1}
 \| V^n_{h,1} + V^n_{h,2} \|_{H^s((S_0)_{\delta} \setminus (S_0)_{h})} \: \le \: C \, h^{\frac{1}{3}-s} \, \|(U^n,U^n_S)\|_{H^1((S_0)_{\delta} \setminus (S_0)_{h}) \times {\cal R}}. 
\end{equation} 
Finally, we build some  $W^n_h = \nabla Y^n_h$ where $Y^n_h$ is the unique solution of :
$$
\left\{
\begin{array}{rcll}
\Delta Y^n_h &=& 0\, & \text{in $(S_0)_{\delta} \setminus (S_0)_{h}$}\,, \\
\partial_{z} Y^n_h  &=& (U^n_S - U^n)\cdot e_z\,, & \text{on $\partial(S_0)_{h}$}\,, \\
\partial_{z} Y^n &=& 0\,, & \text{on $\partial(S_0)_{\delta}\,,$}
\end{array}
\right. 
\quad
\text{ such that }
\int_{(S_0)_{\delta} \setminus (S_0)_{h}} Y^n_h = 0 \,.
$$ 
we recall that $\nu = e_z$ on $\partial (S_0)_h.$
By standard elliptic regularity results, there exists a constant $C$ independent of $h$ such that :
\begin{eqnarray*}
\|W^n_h\|_{L^2((S_0)_{\delta} \setminus (S_0)_{h})}  &\le & \|Y^n_h\|_{H^1((S_0)_{\delta} \setminus (S_0)_{h})} \leq C\|(U^n_S - U^n)\cdot e_z\|_{H^{-1/2}(\partial(S_0)_{h})}\,,\\
\|W^n_h\|_{H^1((S_0)_{\delta} \setminus (S_0)_{h})}  &\le & \|Y^n_h\|_{H^2((S_0)_{\delta} \setminus (S_0)_{h})} \leq C\|(U^n_S - U^n)\cdot e_z\|_{H^{1/2}(\partial(S_0)_{h})}\,. 
\end{eqnarray*}
By interpolation, we get 
\begin{equation} \label{interpol2}
\|W^n_h\|_{H^{1/2}((S_0)_{\delta} \setminus (S_0)_{h})}  \: \le \:   C \|(U^n_S - U^n)\cdot e_z\|_{L^2(\partial(S_0)_{h})} \\
 \end{equation}
Now, we write
\begin{multline}  \label{interpol3}
 \|(U^n_S - U^n) \cdot e_z\|_{L^2(\partial (S_0)_{h})}  \\
 \leq C \left( h^{\frac{1}{2}} \|\nabla (U^n_S-U^n)\|_{L^2((S_0)_{\delta} \setminus (S_0)_{h})} \: + \:  
 \|(U^n_S - U^n) \cdot e_z\|_{L^2(\partial S_0)} \right).
\end{multline}
Eventually, we set $\displaystyle V^n_h \: := \:  V^n_{h,1} + V^n_{h,2} - W^n_h$. We stress that the normal component of $\displaystyle V^n_h$ is continuous across $\displaystyle \pa (S_0)_h$ and  $\displaystyle \pa (S_0)_\delta$. Hence, for any $s < \frac{1}{2}$, the $H^s$ norm of $V^n_h$ over the whole domain is controlled by the sum of the $H^s$ norms over $\displaystyle (S_0)_h$, $\displaystyle (S_0)_{\delta} \setminus (S_0)_{h}$ and $\displaystyle \Omega^n\setminus (S_0)_\delta$, where $\Omega^n$ is a shorthand for $\phi^n_{0,t}(\Omega)$. It follows from this remark and  the previous inequalities that: for all $s < 1/2$ 
\begin{equation}
\| V^n_h \|_{H^s(\Omega^n)} \: \le \: C \, \left( (1 + \:   h^{\frac{2-5s}{6}}) \left(\,\|U^n \|_{H^1(F^n)} + \| U^n_S \|_{{\cal R}} \right) +  \|(U^n_S - U^n) \cdot e_z\|_{L^2(\partial S_0)}  \right). 
\end{equation}
 Also, one has 
\begin{multline}
\| V^n_h - U^n \|_{L^2(\Omega^n\setminus (S_0)_h)} \: \\ \le \: C \, \left( h^{1/3} \left(\,\|U^n \|_{H^1(F^n)} + \| U^n_S \|_{{\cal R}} \right) 
+  \|(U^n_S - U^n)  \cdot e_z\|_{L^2(\partial S_0)}  \right). 
\end{multline}

\medskip
Back to the moving domain, and accouting for time dependence, we obtain \eqref{u-vh}.

\bigskip
The second step in the treatment of the nonlinear terms is a control of the Hausdorff distance between $S^n(t)$ and $S(t')$ for close times $t,t' \in [0,T]$. 
This is the purpose of 
\begin{lemma} \label{lemmahausdorff}
Let $h > 0$. 
\begin{description} 
\item[i)] There exists $n_0 \ge 0$ such that for all $ n \ge n_0$, 
$$ S^n(t) \: \subset \:  (S(t))_{h/4}\: \subset \: (S^n(t))_{h/2} \quad  \forall\, t \in [0,T]. $$
\item[ii)] There exists $\eta > 0$ such that for all $t_0 \in [0,T]$,  for all $t \in [t_0 - \eta, t_0 + \eta] \cap [0,T]$ 
$$ (S(t))_{h/2} \: \subset \:   (S(t_0))_{h} \subset (S(t))_{2h}.  $$
\end{description}
\end{lemma}
Note that condition \eqref{deltafar} and point i) of the lemma (applied with $h = \delta$) imply that 
\begin{equation} \label{deltafar2}
 \mbox{dist}(S(t), \pa \Omega) \ge \frac{3}{2} \delta, \quad \mbox{ for } \: t \in [0,T], \quad  \mbox{ for some fixed } \: \delta > 0. 
\end{equation}

{\em Proof of the lemma}. We first treat  {\bf i)}, focusing on the first inclusion  (the other one is proved in the same way). To this end, we recall that the associated sequence of characteristic functions 
$\chi_S^n$ converges to $\chi_S$ in $C([0,T] ; L^1(\Omega)).$ This implies that 
$$
\sup_{t \in [0,T]} |S^n(t) \: \triangle \: S(t)| = \sup_{t \in [0,T]} \|\chi_S^n(t,\cdot) - \chi_S(t,\cdot)\|_{L^1(\Omega)} \to 0 \text{ when $n \to \infty$}\,,
$$
where we denoted $\triangle$ the symmetric difference of subsets of $\mathbb R^3.$ Let  us now take $h >0$ and assume {\it a contrario} that 
there exists a sequence of times $t_k \in [0,T]$ and of integers $n_k$ going to infinity such that 
$$
S^{n_k}(t_k) \setminus (S(t_k))_{h/4} \neq \emptyset\,.
$$  
As $S^{n_k}(t_k)$ is isometric to $S_0$, which satisfies: 
$$
\exists \, r  >0 \text{ s.t.  for all}   \, x \in S_0 \text{ there exists a euclidean ball $B$ with radius $r$ satisfying $x \in B \subset S_0$}
$$
 there exists for all $k$ a ball $B'_k$ with radius $r' = \min (r,h/16)$ such that 
$$
B'_k \subset S^{n_k}(t_k) \setminus S(t_k) \,,
$$
so that 
$$
\sup_{t \in [0,T]} |S^{n_k}(t) \: \triangle \: S(t)| \geq \dfrac{4\pi |r'|^3}{3}\,,
$$
which yields a contradiction. Consequently, there exists $n_0$ such that, for all $n \geq n_0,$
$$ S^n(t) \: \subset \:  (S(t))_{h/4}\,, \quad \forall\, t \in [0,T] .$$

The second item  {\bf ii)} is  obtained in the same way. Let $h>0$ and assume for instance that the first inclusion does not hold. 
Arguing as previously, we are  able to construct two sequences $(t_0^k)$ and $(t^k)$ converging both 
to $t_0 \in [0,T]$ and such that $S(t^k) \setminus S(t^k_0) $ contains a ball of fixed radius. Once again,
this contradicts the continuity in $L^1(\Omega)$ of $\chi_S$ at $t_0.$

\bigskip

Thanks to the previous lemmas, we can  conclude the proof of Proposition \ref{prop_strongcvgce}, following very closely \cite{SanMartin&al02}. At first, very minor adaptation of  the  proof of \cite[Proposition 7.1]{SanMartin&al02} yields:    {\em for $s \in (0,1/3),$ there exists $h_0$, such that, for all $h \in (0,h_0)$:
\begin{equation} \label{eq_rhonun}
\lim_{n\to \infty} \int_0^{T} \int_{\Omega} \rho^n u^n \cdot P^{s}[(S(t))_h] u^n =   \int_0^{T} \int_{\Omega} \rho u \cdot P^{s}[(S(t))_h] u.
\end{equation}}
We remind that the main idea behind this limit is the following:  thanks to Lemma \ref{lemmahausdorff}, for any  field $\xi$, the projected field $P^{s}[(S(t))_h](\xi)$ is rigid in a neighborhood of $S^n(t)$ for $n$ large enough. Hence, if one uses   $P^{s}[(S(t))_h](\xi)$  as a test function in the momentum equation, the boundary term at $\pa S^n(t)$ and the penalization term vanish: roughly, one recovers  a uniform bound on $\pa_t P^{s}[(S(t))_h] (\rho^n u^n)$  in a Sobolev space of negative index, and from there compactness. 
For all details, see \cite[Proposition 7.1]{SanMartin&al02}.

\medskip
Then, one establishes that 
\begin{equation} \label{eq_normL2}
\lim_{n\to \infty} \int_{0}^T \int_{\Omega} \rho^n |u^n|^2 = \int_0^{T} \int_{\Omega} \rho |u|^2,
\end{equation}
$\biggl( \rho^n = \rho_F (1 - \chi^n_S) + \rho_S \chi^n_S, \quad \rho = \rho_F (1 - \chi_S) + \rho_S \chi_S\biggr).$ The idea is to write 
\begin{align*} 
&  \int_{0}^T \int_{\Omega} \rho^n |u^n|^2 - \int_0^{T} \int_{\Omega} \rho |u|^2 \\
& = \left( \int_0^{T} \int_{\Omega} \rho^n u^n \cdot P^{s}[(S(t))_h](u^n) - \int_0^{T} \int_{\Omega} \rho u \cdot P^{s}[(S(t))_h](u) \right)  \\
& \quad + \int_0^T \int_\Omega \rho^n u^n \cdot (u^n - P^s[S(t)_h] u^n) dt \: + \: \int_0^T \int_\Omega \rho u \cdot (P^s[S(t)_h] u - u) dt
\end{align*}
The first term  at the r.h.s. is controlled using \eqref{eq_rhonun}, whereas the last two are treated thanks to Lemma \ref{prop_strongcvgce}: note that $(S(t))_h  \subset (S^n(t))_{2h}$ for $n$ large enough  by Lemma \ref{lemmahausdorff},  so that
\begin{align*}
\int_0^T \|u^n(t,\cdot) - P^s[(S(t))_h] \,  u^n(t,\cdot)\|^2_{H^s(\Omega)}\: & \le \: \int_0^T \|u^n(t,\cdot) - P^s[(S^n(t))_{2h}] \,  u^n(t,\cdot)\|^2_{H^s(\Omega)} \\
&  \: \leq  \: C \left(h^\eps \, + \, n^{-\eps} \right). 
 \end{align*}

\medskip
The final step of the proof consists in showing that 
$$ \int_0^T \int_\Omega \rho |u^n|^2 \: \rightarrow \: \int_0^T \int_\Omega \rho |u|^2 $$
which yields the strong compactness of $u^n$ ($\rho$ having positive lower and upper bounds).
The idea here is to write 
$$ \left| \int_0^T \int_\Omega \rho ( |u^n|^2 - |u|^2 ) \right| \: \le \:  
\left|  \int_0^T \int_\Omega \left( \rho^n |u^n|^2 - \rho |u|^2 \right) \right| \: + \: \left|    \int_0^T \int_\Omega (\rho^n - \rho) |u^n|^2 \right|. $$
The first term at the r.h.s. goes to zero by  \eqref{eq_normL2}. For the second one, we use that $\rho^n\rightarrow \rho$ strongly in $C([9,T];L^p(\Omega))$ for all finite $p$ and that $|u^n|^2$ is uniformly bounded in $L^{p'}$ for some $p' > 1$, thanks to the uniform $H^s$ bound on $u^n$. 
Again, we refer to \cite{SanMartin&al02} for all details.

\subsection{Convergence in the momentum equation: nonlinear terms}

 \medskip
 Thanks to the strong convergence of {Proposition \ref{prop_strongcvgce}}, we are now able to split ${{\rm conv}^n}$ in a suitable way. 
Let us first remind that $v^n$ is identically equal to $u^n_S$ inside $S^n$, whereas $\varphi^n$ is identically equal to $\varphi_F$ outside $S^n$. This allows us to decompose the convection term as follows:
\begin{align*}
{\rm conv}^n & =  \int_0^T \int_\Omega \rho_F (1- \chi^n_S) \, u^n_F \cdot  \pa_t \varphi_F  \\
& + \int_0^T \int_\Omega  \rho_F (1- \chi^n_S)  v^n \otimes u^n_F : \na \varphi_F  + \int_0^T \rho_S \chi^n_S (\pa_t + u^n_S \cdot \na) \varphi^n \cdot u^n \: := \: I^n_1 + I^n_2 + I^n_3 
\end{align*} 
The convergence of $I_1^n$ is clear: 
$$ I_1^n \rightarrow  \int_0^T \int_\Omega \rho_F (1- \chi_S) \, u_F \cdot  \pa_t \varphi_F. $$
The convergence of $I_3^n$ follows from the fourth item in Proposition \ref{approximatetest}, which clearly implies that 
$$ I^n_3 = \int_0^T \int_\Omega  \rho_S \chi^n_S (\pa_t + u^n_S \cdot \na) \varphi_S \cdot u^n + o(1).  $$ 
Using the strong convergence of $\chi^n_S u^n$ to $\chi_S u_S$ in $L^{2}((0,T) \times \Omega)$,  it is then easily shown that 
$$ I^n_3  =  \int_0^T \int_\Omega  \rho_S \chi_S \pa_t  \varphi_S \cdot u_S + \int_0^T \int_\Omega  \rho_S \chi^n_S u^n_S \cdot \na \varphi_S \cdot u^n_S \: + \:  o(1).$$ 
Now, we write the second term at the r.h.s. as 
\begin{align*}
 \int_0^T \int_\Omega  \rho_S \chi^n_S u^n_S \cdot \na \varphi_S \cdot u^n_S & \: = \:    \int_0^T \int_\Omega  \rho_S \chi^n_S u^n_S \otimes u^n_S : \na \varphi_S \\
  & \: = \:    \int_0^T \int_\Omega  \rho_S \chi^n_S u^n_S \otimes u^n_S : D(\varphi_S) \: = \:  0 
  \end{align*}
  as $\varphi_S$ is a rigid vector field. 
 
 \medskip
 It remains to study  $I^n_2$. We know from paragraph \ref{aprioriboundsv}  that 
 $$ (1-\chi^n_S) (v^n - u^n)  =  (1-\chi^n_S) (v^n - u^n_F) \rightarrow 0 \: \mbox{ in } \:  L^2(0,T; L^p(\Omega)), \quad \forall p \le 6. $$
 It follows that 
\begin{eqnarray*}  
I^n_2 & \: = \: &     \int_0^T \int_\Omega  \rho_F (1- \chi^n_S)  u^n_F \otimes u^n_F : \na \varphi_F + o(1).  \\
	 & \: = \: &    \int_0^T \int_\Omega  \rho_F (1- \chi^n_S)  u^n \otimes u^n : \na \varphi_F + o(1).
\end{eqnarray*}
In this last identity we collect the strong convergences of $u^n$ to $u$ in $L^2((0,T) \times \Omega)$ and of
$\chi^n_S$ to $\chi$ in $\mathcal{C}([0,T];L^{15}(\Omega)),$  together with the
uniform regularity of $(u^n,u)$ in $L^2(0,T;H^{1/5}(\Omega))$ (see Lemma \ref{lemvhbar}), which yields that $(u_n,u)$
are uniformly bounded in $L^2(0,T ; L^{30/13}(\Omega)).$  We obtain then :
\begin{eqnarray*}
\lim_{n\to\infty} \int_0^T \int_\Omega  \rho_F (1- \chi^n_S)  u^n \otimes u^n : \na \varphi_F &=& \int_0^T \int_\Omega  \rho_F (1- \chi_S)  u\otimes u : \na \varphi_F\\
					&=& \int_0^T \int_\Omega  \rho_F (1- \chi_S)  u_F\otimes u_F : \na \varphi_F
\end{eqnarray*}
This concludes our proof.

\subsection{Energy inequality and extension to collision time.}
We pass to the weak limit in \eqref{energyboundn} and prove that the solution $(\rho,u)$ satisfies the further energy estimate \eqref{mainenergybound}.
First, we note that \eqref{energyboundn} implies 
\begin{align} & \| \sqrt{\rho^n} u^n(t, \cdot) \|^2_{L^2(\Omega)} \: + \:  \int_0^t \int_\Omega  2 \mu^n |D(u^n)|^2 \: + \:  \frac{1}{2\beta_\Omega}\int_0^t \int_{\pa \Omega}  | u^n \times \nu|^2 \notag  \\
& \: + \:  \frac{1}{2\beta_S} \int_0^t  \int_{\pa S^N(t)}  |(u^n - P^n_S u^n) \times \nu|^2  \le \: \int_0^t  \int_{\Omega} \rho (-g) \cdot u^n  \: + \:  \int_{\Omega} \rho_0\, |u_0|^2 \notag
\end{align}
for all $n.$  As we have convergence of $\sqrt{\rho^n} u^n$ in $L^{2-\varepsilon}((0,T) \times \Omega)$  we can pass to the weak limit in 
this inequality for almost all $t \in [0,T].$ As $S(t)$ remains far from $\partial \Omega$ we treat boundary terms in a similar way as in paragraph \ref{sec_nrjin}.
The only term which requires a new treatment is :
$$
\int_0^t \int_\Omega  2 \mu^n |D(u^n)|^2.
$$
For this term, we note that, because of Lemma \ref{lemmahausdorff}, there holds  for arbitrary $h >0$ and  $n$ sufficiently large :
$$
\int_0^t \int_{\Omega \setminus (S(t))_h}  2 \mu_F |D(u^n)|^2 \leq \int_0^t \int_\Omega  2 \mu^n |D(u^n)|^2.
$$
If we let $n$ go to infinity, and then $h$   go to $0$, we obtain :
$$
\int_{0}^t \int_{F(t)}  2 \mu_F |D(u_F)|^2 \leq \liminf \int_0^t \int_\Omega  2 \mu^n |D(u^n)|^2,
$$
for almost all $t \in [0,T].$  Hence, passing to the limit in \eqref{energyboundn} yields \eqref{mainenergybound}.

\medskip

Our solutions are limited in time to avoid collision. Namely, the only shortcoming of our construction is that it requires the distance
between $S(t)$ and $\partial \Omega$ to be larger than a fixed positive distance $\delta$ through time. However, as long as we are 
given an initial data $u_0 \in L^2(\Omega)$ and an initial position $S_0$ such that $S_0  \Subset \Omega,$ we are able to construct 
a small time $T$ depending only on the inital position of $S_0$ in $\Omega$ and the $L^2$ norm of $u_0$ such that the solution exists
and satisfies this  property on $[0,T].$  As our solutions satisfy also energy estimate \eqref{mainenergybound} we might reproduce the arguments of \cite[Lemma 2.2]{Feireisl03} to concatenate  solutions in time and prove existence of at least one weak solution until collision time.

\appendix

\section{Weak/Strong convergence and isometries} \label{app_lemmawn}
In this appendix, we study the influence of isometric transformations on weak and strong convergence of sequences.
First, we prove :
\begin{lemma} \label{lemmaw}
Let $\phi \in C([0,T] ; {\mbox Isom}(\R^3)).$ Given $w: (0,T) \times \mathbb R^3 \to \mathbb R^3$ 
we define : 
\begin{equation} \label{eq_wW}
w(t,\phi(t,y)) := d\phi_t|_y W(t,y)\,, \quad \forall \, (t,y) \in (0,T) \times \mathbb R^3. 
\end{equation}
Then
\begin{itemize}
\item If $w \in L^2(0,T ; H^1(\mathbb R^3))$ then $W \in L^2(0,T;H^1(\mathbb R^3)).$
\item If $w \in C([0,T] ; H^1(\mathbb R^3))$ then $W \in C([0,T] ; H^1(\mathbb R^3)).$ 
\item The same assertions hold true replacing $H^1(\mathbb R^3)$ by $H^1_{loc}(\mathbb R^3).$
\end{itemize}
\end{lemma}
The proof of this lemma is based on the fact that formula \eqref{eq_wW} for fixed $t$ defines a unitary
transformation of $H^1(\mathbb R^3).$ The details are left to the reader. Second, we obtain :
 \begin{lemma} \label{lemmawn} 
Let $\phi^N : [0,T] \times \mathbb R^3$ such that $\phi^N(t,\cdot) \in {\mbox Isom}(\mathbb R^3)$ for all $t \in [0,T].$
We assume that $\phi^N$  converges to $\phi$ in $C([0,T];C^{\infty}_{loc}(\mathbb R^3)).$ 
Given a sequence $(w^N): (0,T) \times \mathbb R^3 \to \mathbb R^3$ we define : 
$$
w^N(t,\phi^N(t,y)) := d\phi^N_t\vert_y W^N(t,y).
$$ 
Then, with obvious notations:
\begin{itemize}
\item If $(w^N)$ converges to $w$ strongly (resp. weakly) in $L^2(0,T ; H^1(\mathbb R^3))$ then $(W^N)$
converges to $W$ strongly (resp. weakly) in $L^2(0,T;H^1(\mathbb R^3)).$
\item  If $(w^N)$ converges to $w$  in $C([0,T] ; H^1(\mathbb R^3))$ then $(W^N)$
converges to $W$ in $C([0,T];H^1(\mathbb R^3)).$
\item The same assertions hold true replacing $H^1(\mathbb R^3)$ by $H^1_{loc}(\mathbb R^3)$.
\end{itemize}
\end{lemma}

\medskip
\noindent
{\em Remark.} We point out that $w^N$ and $W^N$ satisfy symmetric relations: 
$$ w^N\left(t, \phi^N_{t}(y)\right)\: = \:  d\phi^N_{t}\vert_{y} \, W^N(t,y) \quad  \Leftrightarrow  \quad 
W^N\left(t, [\phi^N_{t}]^{-1}(x)\right)\: = \:  d[\phi^N_{t}]^{-1}\vert_{x} \, w^N(t,x), $$
so that fields $W^n$ and $w^n$, resp. $W$ and $w$ can be switched in this lemma. 

\medskip

\noindent{\em Proof of Lemma \ref{lemmawn}}. 
 We first remind that $\phi^N_{t}$ is an affine isometry, so that  (for all $N,t$)
\begin{equation} \label{isometry}
\big|d\phi^N_{t}\vert_{y} \, x \big| = |x|\,,  \quad \big|[d\phi^N_{t}\vert_{y}]^{-1} \, M \, d\phi^N_{t}\vert_{y}\big| = |M|, \quad \forall (x,y) \in \R^3 \times \mathbb R^3, \quad \forall M \in M_3(\R). 
\end{equation}
The same relations hold for $\phi$ instead of $\phi^N$.  

\medskip

{\em Strong convergence.}
We  focus on convergence in $C([0,T] ; H^1(\mathbb R^3))$, the  strong convergence in $L^2 H^1$ being treated in the same way. 
First, we note that our previous lemma yields: 
$$W^N \in C([0,T];H^1(\mathbb R^3)) \:  \mbox{ for any } N, \quad   W \in C([0,T];H^1(\mathbb R^3)).$$
 Then, we write
\begin{equation*}
\| W^N - W \|_{C([0,T] ;  L^2(\R^3))}   \: \le \: \sup_{t} I_1^N(t) \: + \sup_{t} \: I_2^N(t)    \: + \sup_{t}  \: I_3^N(t)
\end{equation*}
where
$$ |I_1^N(t)|^2  :=  \int_{\R^3}  \left|  W^N(t,y) - d\phi^N_{t}\circ [d\phi_{t}]^{-1} \vert_y \,  W(t, [\phi_{t}]^{-1} \circ \phi^N_{t}(y)) \right|^2 dy ,  $$
$$ |I_2^N(t)|^2  :=  \int_{\R^3}  \left|   [d\phi^N_{t}]^{-1}\circ d\phi_{t}\vert_y \, W(t, [\phi_{t}]^{-1} \circ \phi^N_{t}(y)) -  W(t,[\phi_{t}]^{-1} \circ \phi^N_{t}(y)) \right|^2 dy, $$
and 
$$ |I_3^N(t)|^2  :=  \int_{\R^3} \left|   W(t,[\phi_{t}]^{-1} \circ \phi^N_{t}(y)) - W(t,y) \right|^2 dy . $$
 Using \eqref{isometry}, we have easily 
 \begin{align*} |I_1^N(t)|^2 & \: = \: \int_{\R^3} \left| w^N(t,\phi^N_{t}(y)) - w(t,\phi^N_{t}(y)) \right|^2 dy dt \\
& \: = \:  \int_{\R^3}  \left| w^N(t,x) - w(t,x) \right|^2 dx,
\end{align*}
which tends uniformly to $0$ when $N\to \infty$ by assumption. 
We then get
\begin{align*} |I_2^N(t)|^2   \: & \le \: \sup_{t,y}  \left|  d\phi^N_{0,t}\circ d\phi_{t,0} - I_d \right|^2 \,  \int_{\R^3}  \left|  W(t, [\phi_{t}]^{-1} \circ \phi^N_{t}(y)) \right|^2 dy  \\
  \: & \le \: \sup_{t,y}  \left|  d\phi^N_{0,t}\circ d\phi_{t,0} - I_d \right|^2 \| W \|_{C([0,T];H^1(\mathbb R^3))} \rightarrow 0. 
  \end{align*}
Finally, the continuity of $W$ with values in $L^2(\mathbb R^3)$ implies that:
$$
\int_{|y| \ge A} | W(t,y)|^2 dy 
$$ 
can be made arbitrary small uniformly in time, taking $A$ sufficiently large. So, we apply the local convergence
of $\phi^N$ to $\phi$ to obtain that, for $N$ sufficiently large, there holds :
\begin{align*}
 |I_3^N(t)| \: & \le \:  \left( \int_{|y| \ge A} | W(t,[\phi_{t}]^{-1} \circ \phi^N_{t}(y)) |^2 dy\right)^{1/2} \: + \: 
   \left(   \int_{|y| \ge A} | W(t,y) |^2 dy \right)^{1/2}   \\
   & + \:  \left( \int_{|y|  <  A}   \left|   W(t,[\phi_{t}]^{-1} \circ \phi^N_{t}(y)) - W(t,y) \right|^2 dy  \right)^{1/2}  \\
&  \le \: 2     \left(  \int_0^T \int_{|y| \ge A/2} | W(t,(y)) |^2 \right)^{1/2} + \:  \left(  \int_{|y| < A}   \left|   W(t,[\phi_{t}]^{-1} \circ \phi^N_{t}(y)) - W(t,y) \right|^2 dy \right)^{1/2}  
\end{align*}
The first term at the r.h.s. is independent of $N$ and goes to zero as $A$ goes to infinity. Moreover, for fixed $A$,  $[\phi_{t}]^{-1} \circ \phi^N_{t}(y)$ converges to $y$ uniformly in $[0,T] \times \{ |y| \le A \}$. Hence, for fixed $A$, the second term at the r.h.s. converges to zero as $N$ goes to infinity (continuity of translations in $L^2$) uniformly in $t$. We conclude that $I_3^N$ goes to zero, so that $W^N$ converges to $W$ in $C([0,T] ; L^2(\R^3))$. The convergence of  $\na W^N$ to $\na W$ follows the same lines, which yields the result. 

\medskip
{\em Weak convergence}. Again, we only prove the convergence on $\mathbb R^3.$ The convergence in $H^1_{loc}(\mathbb R^3)$ is similar.
We assume here that 
 $(w^N)$ converges to $w$ in $L^2(0,T; H^1(\R^3))$ weak.  Given $\chi \in C^{\infty}_c((0,T) \times \mathbb R^3)$ 
there holds :
\begin{eqnarray*}
\int_{0}^T \int_{\mathbb R^3} W^N(t,y) \cdot \chi(t,y) \: dt \: dy 
&=& \int_{0}^T \int_{\mathbb R^3}  w^N(t,\phi_{t}^N(y)) \cdot  \, d\phi^N_{t,0}\vert_{y} \,  \chi(t,y) \: dt \: dy  \\
&=& \int_{0}^T \int_{\mathbb R^3}   w^N(t,x) \cdot  \left(d[\phi^N_{t}]^{-1}\vert_x \right)^{-1} \chi(t,[\phi^N_{t}]^{-1}(x)) \: dt \: dy\,, 
\end{eqnarray*}
where we applied again that $d\phi^N_{t}\vert_{y}$ is a linear isometry. 
Because of the strong convergence of $\phi^N$ in  $C([0,T]; C^1_{loc}(\mathbb R^3))$
there holds :
$$
\left(d[\phi^N_{t}]^{-1}\vert_x \right)^{-1} \chi(t,[\phi^N_{t}]^{-1}(x)) \to \left(d[\phi_{t}]^{-1}\vert_x \right)^{-1} \chi(t,[\phi_{t}]^{-1}(x))\text{ strongly in $L^{2}((0,T) \times \mathbb R^3)$}
$$
so that  with the weak convergence of $w^N$ we obtain :
$$
\int_{0}^T \int_{\mathbb R^3} W^N(t,y) \cdot \chi(t,y) \: dt \: dy  \to \int_{0}^T \int_{\mathbb R^3} W(t,y) \cdot \chi(t,y) \: dt \: dy\,.
$$
Similar arguments yield also that:
$$
\int_{0}^T \int_{\mathbb R^3} \nabla W^N(t,y) :  \nabla \chi(t,y) \: dt \: dy  \to \int_{0}^T \int_{\mathbb R^3} \nabla W(t,y) : \nabla \chi(t,y) \: dt \: dy\,.
$$
which ends the proof.


\begin{thebibliography}{99}

\bibitem{Barnocky}
G.~Barnocky and R.~H. Davis.
\newblock The influence of pressure-dependent density and viscosity on the
  elastohydrodynamic collision and rebound of two spheres.
\newblock {\em J. Fluid Mech.}, 209:501--519, 1989.


\bibitem{Bocquet:2007}
{C. Ybert, C. Barentin, C. Cottin-Bizonne, P. Joseph, and L. Bocquet}
\newblock Achieving large slip with superhydrophobic surfaces: Scaling laws
for generic geometries 
\newblock {\em Physics of fluids} 19, 123601 (2007).

\bibitem{BostCottetMaitre10}
C.~Bost, G.-H. Cottet, and E.~Maitre.
\newblock Convergence analysis of a penalization method for the
  three-dimensional motion of a rigid body in an incompressible viscous fluid.
\newblock {\em SIAM J. Numer. Anal.}, 48(4):1313--1337, 2010.


\bibitem{Brenner&Cox63}
H.~Brenner et R.~G. Cox.
\newblock The resistance to a particle of arbitrary shape in translational
  motion at small {R}eynolds numbers.
\newblock {\em J. Fluid Mech.} {\bf 17} (1963),  561--595.


\bibitem{Cooley&Oneill69}
M.~Cooley and M.~O'Neill.
\newblock On the slow motion generated in a viscous fluid by the approach of a
  sphere to a plane wall or stationary sphere.
\newblock {\em Mathematika}, 16:37--49, 1969.




\bibitem{Davis3}
R.~H. Davis, J.~Serayssol, and E.~Hinch.
\newblock The elastohydrodynamic collision of two spheres.
\newblock {\em J. Fluid Mech.}, 163:045302, 2006.

\bibitem{Davis1}
R.~H. Davis, Y.~Zhao, K.~P. Galvin, and H.~J. Wilson.
\newblock Solid-solid contacts due to surface roughness and their effects on
  suspension behaviour.
\newblock {\em R. Soc. Lond. Philos. Trans. Ser. A Math. Phys. Eng. Sci.},
  361(1806):871--894, 2003.



\bibitem{Desjardins&Esteban99}
B.~Desjardins et M.~J. Esteban, 
\newblock Existence of weak solutions for the motion of rigid bodies in a
  viscous fluid.
\newblock {\em Arch. Ration. Mech. Anal.} {\bf 146} (1999), 59--71.



\bibitem{DiPernaLions89}
R.~J. DiPerna et P.-L. Lions.
\newblock Ordinary differential equations, transport theory and {S}obolev
  spaces.
\newblock {\em Invent. Math.} {\bf  98} (1989), 511--547.


 
\bibitem{Feireisl03}
E. Feireisl.
\newblock On the motion of rigid bodies in a viscous incompressible fluid.
\newblock {\em J. Evol. Equ.}, {\bf 3} (2003), 419--441.
 
\bibitem{Galdibooknew}
G.~P. Galdi.
\newblock {\em An introduction to the mathematical theory of the
  {N}avier-{S}tokes equations}.
\newblock Springer Monographs in Mathematics. Springer, New York, second
  edition, 2011.
\newblock Steady-state problems.


\bibitem{GeHi}
M. Hillairet et D.   G\'erard-Varet. 
\newblock Regularity issues in fluid structure interaction
\newblock {\em Arch. Rat. Mech. Anal.}, {\bf 195} (2010), 375--407. 

\bibitem{GeHi2}
D. G\'erard-Varet and M. Hillairet.
\newblock Computation of the drag force on a sphere close to a wall. The roughness issue.
\newblock {\em M2AN}, to appear


\bibitem{GlassSueurpp}
O.~Glass and F.~Sueur.
\newblock Uniqueness results for weak solutions of two-dimensional fluid-solid
  systems.
\newblock arXiv:1203.2894v1, March 2012.

\bibitem{henrot} A. Henrot and M. Pierre 
\newblock Variation et optimisation de formes. Une analyse g\'eom\'etrique (French) [Shape variation and optimization. A geometric analysis]
\newblock Math\'ematiques $\&$ Applications 48, Springer, Berlin, 2005.



\bibitem{Hillairet07}
M.~Hillairet.
 \newblock Lack of collision between solid bodies in a 2{D} incompressible
  viscous flow.
\newblock {\em Comm. Partial Differential Equations}, {\bf 32} (2007), 1345--1371.


\bibitem{HiTa}
M. Hillairet et T.  Takahashi.
\newblock Collisions in 3D fluid structure interaction problems
\newblock {\em SIAM J. Math. Anal.}, {\bf 40} (2009), 2451--2477. 



\bibitem{Hocking}
L.~Hocking.
\newblock The effect of slip on the motion of a sphere close to a wall and of
  two adjacent sheres.
\newblock {\em J. Eng. Mech.}, 7:207--221, 1973.


\bibitem{Lauga:2007}
{E. Lauga, M.P. Brenner and H.A. Stone}
\newblock Microfluidics: The no-slip boundary condition
\newblock {\em Handbook of Experimental Fluid Dynamics},C. Tropea,
A. Yarin, J. F. Foss (Eds.), Springer, 2007. 



\bibitem{PLLvol1}
P.-L. Lions.
\newblock {\em Mathematical topics in fluid mechanics. {V}ol. 1}, volume~3 of
  {\em Oxford Lecture Series in Mathematics and its Applications}.
\newblock The Clarendon Press Oxford University Press, New York, 1996.
\newblock Incompressible models, Oxford Science Publications.
 

\bibitem{MaPu}
C. Marchioro and R. Pulvirenti.
\newblock Mathematical Theory of Incompressible Nonviscous fluids
\newblock {\em Springer-Verlag, New York}, Applied Mathematical Sciences vol. 96, 1994. 

\bibitem{Nitsche:1981}
{\sc Nitsche, J. A.}
\newblock On Korn's second inequality.
\newblock {\em R.A.I.R.O. Analyse num\'erique/Numerical Analysis} 15 (1981), 237--348.



\bibitem{SanMartin&al02}
J.A. San~Mart{\'{\i}}n, V.~Starovoitov, et M.~Tucsnak.
\newblock Global weak solutions for the two-dimensional motion of several rigid
  bodies in an incompressible viscous fluid.
\newblock {\em Arch. Ration. Mech. Anal.}, {\bf 161} (2002), 113--147, 2002.

\bibitem{Smart}
J.~Smart and D.~Leighton.
\newblock Measurement of the hydrodynamic surface roughness of noncolloidal
  spheres.
\newblock {\em Phys. Fluids}, 1:52--60, 1989.


\bibitem{Starovoitov03}
V.~N. Starovoitov. 
\newblock Behavior of a rigid body in an incompressible viscous fluid near a
  boundary.
\newblock In {\em Free boundary problems (Trento, 2002)}, volume 147 of {\em
  Internat. Ser. Numer. Math.},   313--327. Birkh\"auser, Basel, 2004.

\bibitem{Starovoitov03bis}
V.~N. Starovoitov.
\newblock Nonuniqueness of a solution to the problem on motion of a rigid body
  in a viscous incompressible fluid.
\newblock {\em J. Math. Sci.}, 130(4):4893--4898, 2005.

\bibitem{stokes}  G. Stokes.
 On the Theories of Internal Friction of 
Fluids in Motion and of the Equilibrium and Motion of Elastic Solids, 
{\em Trans. Camb. Phil. Soc.}, {\bf 8} (1845),    287-319.

\bibitem{Takahashi03bis}
T.~Takahashi.
\newblock Analysis of strong solutions for the equations modeling the motion of
  a rigid-fluid system in a bounded domain.
\newblock {\em Adv. Differential Equations}, {\bf 8} (2003):1499--1532.

\bibitem{Vinogradova1}
O.~Vinogradova and G.~Yakubov.
\newblock Surface roughness and hydrodynamic boundary conditions.
\newblock {\em Phys. Rev. E}, 73:479--487, 1986.

\end{thebibliography}
\end{document}